\documentclass[12pt]{article}
\usepackage[top=4cm,bottom=3cm,left=4cm,right=3cm]{geometry}
\usepackage[square,numbers]{natbib}
\usepackage{amssymb}
\usepackage{amsthm}

\usepackage{tabularx,environ}
\usepackage{soul}
\usepackage{caption}
\usepackage{subcaption}
\usepackage{tikz}
\usepackage{amsmath}
\usepackage{multirow}
\usepackage{color}

\usepackage{rotating}
\usepackage{fixmath}
\usepackage[noend]{algpseudocode}
\usepackage[ruled, lined, longend, linesnumbered]{algorithm2e}
\algnewcommand{\IfThenElse}[3]{
  \State \algorithmicif\ #1\ \algorithmicthen\ #2\ \algorithmicelse\ #3}
  \newcommand{\lIfElse}[3]{\lIf{#1}{#2 \textbf{else}~#3}}
\usepackage{lscape}
\usepackage{hyperref}
\usepackage{amssymb}
\usepackage[makeroom]{cancel}
\usepackage[affil-it]{authblk} 
\makeatletter
\newcommand{\problemtitle}[1]{\gdef\@problemtitle{#1}}
\newcommand{\probleminput}[1]{\gdef\@probleminput{#1}}
\newcommand{\problemquestion}[1]{\gdef\@problemquestion{#1}}

\NewEnviron{proble}{
  \problemtitle{}\probleminput{}\problemquestion{}
  \BODY
  \par\addvspace{.5\baselineskip}
  \noindent
  \begin{tabularx}{\textwidth}{@{\hspace{\parindent}} l X c}
    \multicolumn{2}{@{\hspace{\parindent}}l}{\@problemtitle} \\
    \textbf{Input:} & \@probleminput \\
    \textbf{Question:} & \@problemquestion
  \end{tabularx}
  \par\addvspace{.5\baselineskip}
}
\newcommand\notsotiny{\@setfontsize\notsotiny{6.3}{6.3}}
\makeatother

\begin{document}


\title{A multimodal tourist trip planner integrating road and pedestrian networks}

\author[1]{Tommaso Adamo}
\author[2]{Lucio Colizzi}
\author[2]{Giovanni Dimauro}
\author[1]{Gianpaolo Ghiani}

\author[1]{Emanuela Guerriero}

\affil[1]{Dipartimento di Ingegneria dell'Innovazione,  Universit\`{a} del Salento- Italy }
\affil[2]{Dipartimento di Informatica - Universit\`{a} di Bari- Aldo Moro - Italy}

\maketitle
\begin{abstract}
The \textcolor{black}{\textit{Tourist Trip Design Problem}} aims to prescribe a sightseeing plan that maximizes tourist satisfaction \textcolor{black}{while taking} into account a multitude of parameters and constraints, such as \textcolor{black}{the distances among points of interest, the expected duration of each visit, the opening hours of each attraction, the time available daily. In this article we deal with a variant of the problem in which} the mobility environment consists of a pedestrian network and a road network. \textcolor{black}{So, one plan includes a car tour with a number of stops from which pedestrian subtours to attractions (each with its own time windows) depart. We study the problem and develop a method to evaluate the feasibility of solutions in constant time, to speed up the search. This result is used to devise an ad-hoc \textit{iterated local search}}. Experimental results show that our approach can handle realistic instances with up to 3643 points of interest \textcolor{black}{(over a seven day planning horizon)} in few seconds.
\end{abstract}

%


\section{Introduction} \label{intro}
\textcolor{black}{The tourism industry} is one of the fast-growing sectors \textcolor{black}{in} the world. On the wave of digital transformation, \textcolor{black}{this sector} is experiencing a shift from mass tourism to personalized travel. Designing a \textcolor{black}{tailored} tourist trip is a rather complex and time-consuming \textcolor{black}{process}. \textcolor{black}{Therefore}, the \textcolor{black}{use} of expert and intelligent systems \textcolor{black}{can be beneficial}. Such systems typically appear in the form of ICT integrated solutions that perform \textcolor{black}{(usually on a hand-held device)} three main services: recommendation of attractions (Points of Interest, PoIs), route generation and itinerary customization \citep{gavalas2014survey}.
In this research work, we focus on route generation, known in literature as \textcolor{black}{the \textit{Tourist Trip Design Problem}} (TTDP). The objective of the TTDP is to select PoIs that maximize tourist satisfaction, while taking into account a set of parameters (e.g., alternative transport modes, distances among PoIs) and constraints (e.g.\textcolor{black}{, the duration of each visit,} the opening hours of each PoI and \textcolor{black}{the time available daily for sightseeing)}.  In last \textcolor{black}{few} years there has been a flourishing of scholarly work on the TTDP \citep{ruiz2022systematic}. \textcolor{black}{Different variants of TTDP have been studied in the literature, the main classification being made w.r.t. the mobility environment which can be \textit{unimodal} or \textit{multimodal} \citep{ruiz2021tourist}}. 

In this article, we focus on a variant of \textcolor{black}{the} TTDP in which a tourist can move from one PoI to the next one as a pedestrian or as a driver of a vehicle  (like a car or a motorbike).  \textcolor{black}{Under this hypothesis, one plan includes a car tour with a number of stops from which pedestrian subtours to attractions (each with its own time windows) depart.} We refer to this multimodal setting as a \textit{walk-and-drive} mobility environment. \textcolor{black}{Our research work was motivated by a  project aiming to stimulate tourism in the Apulia region (Italy). Unfortunately, the public transportation system is not well developed in this rural area and most attractions can be conveniently reached only by car or scooter, as reported in a recent newspaper article \citep{alongrustyroads}: \textit{(in Apulia) sure, there are trains and local buses, but using them exclusively to cross this varied region is going to take more time than most travellers have.} Our research was also motivated by the need to maintain social distancing in the post-pandemic era \citep{li2020coronavirus}.}

\indent \textcolor{black}{The \textit{walk-and-drive} variant of the TTDP addressed in this article presents several peculiar algorithmic issues that we now describe}. The TTDP is a variant of the \textcolor{black}{\textit{Team Orienteering Problem with Time Windows}} (TOPTW), \textcolor{black}{which is known to be NP-hard \citep{GAVALAS201536}}. \textcolor{black}{We now review the state-of-the-art of modelling approaches, solution methods and planning applications for tourism planning. A systematic review of all the relevant literature has been recently published in \cite{ruiz2022systematic}. The TTDP is a variant of the Vehicle Routing Problem (VRP) with Profits \cite{archetti2014chapter}, a generalization of the classical VRP where the constraint to visit all customers is relaxed. A known profit is associated with each demand node. Given a fixed-size fleet of vehicles, VRP with profits aims to maximize the profit while minimizing the traveling cost. The basic version with only one route is usually presented as Traveling Salesman Problem (TSP) with Profits \cite{feillet2005traveling}. Following the classification introduced in \cite{feillet2005traveling} for the single-vehicle case, we distinguish three main classes. The first class of problems is composed by the Profitable Tour Problems (PTPs) \cite{dell1995prize} where the objective is to maximize the difference between the total collected profit and the traveling cost. The capacitated version of PTP is studied in \cite{archetti2009capacitated}. The second class is formed by Price-Collecting Traveling Salesman Problems (PCTSPs) \cite{balas1989prize} where the objective is to minimize the total cost subject to a constraint on the collected profit. The Price-Collecting VRPs has been introduced in  \cite{tang2006iterated}. Finally, the last class is formed by the \textit{Orienteering Problems} (OPs) \cite{golden1987orienteering} (also called Selective TSPs \cite{laporte1990selective} or Maximum Collection Problems \cite{kataoka1988algorithm}) where the objective is to maximize collected profit subject to a limit on the total travel cost. The \textit{Team Orienteering Problem} (TOP) proposed by \cite{chao1996team} is a special case of VRP with profits; it corresponds to a multi-vehicle extension of OP where a time constraint is imposed on each tour. \\
\indent For the TTDP, the most widely modelling approach is the TOP. Several variants of TOP have been investigated with the aim of obtaining realistic tourist planning. Typically PoIs have to be visited during opening hours, therefore the best known variant is the Team Orienteering Problem with Time-Windows (TOPTW) (\cite{vansteenwegen2009iterated} \cite{boussier2007exact}, \cite{RIGHINI20091191}, \cite{6004465}). In many practical cases, PoIs might have multiple time windows. For example, the tourist attraction is open between 9 am and 14 am and between 3 pm and 7 pm. In \cite{tricoire2010heuristics}, the authors devise a polynomial-time algorithm for checking feasibility of multiple time windows. The size of the problem is reduced in a preprocessing phase if the PoI-based graph satisfies the triangle inequality. 
The model closest to the one proposed in this work is the Multi-Modal TOP with Multiple Time Windows (MM-TOPMTW) \cite{ruiz2022systematic}. Few contributions deal with TTDP in a multimodal mobility environment. Different physical networks and modes of transports are incorporated according to two different models. The former implicitly incorporates multi-modality by considering the public transport. Due to the waiting times at boarding stops, the model is  refereed to as Time-Dependent TOPTW (\cite{zografos2008algorithms}, \cite{garcia2013integrating}, \cite{gavalas2015heuristics}). Other models \textcolor{black}{incorporate}  the choice of transport modes more explicitly, based on availability, preferences and time constraints . In particular in the considered transport modes the tourist either walks or takes a vehicle as passenger, i.e. bus, train, subway, taxi \cite{RUIZMEZA2021107776},\cite{ruiz2021tourist},  \cite{YU20171022}. To the best of our knowledge this is the first contribution introducing the TTDP in a \textit{walk-and-drive} mobility environment. Other variants have been proposed to address realistic instances. Among the others, they include: time dependent profits (\cite{vansteenwegen2019orienteering}, \cite{YU2019488}, \cite{gundling2020time}, \cite{KHODADADIAN2022105794}), score in arcs (\cite{VERBEECK201464}), tourist experiences (\cite{zheng2019using},\cite{RUIZMEZA2021107776},\cite{ruiz2021tourist},\cite{ruiz2021grasp}), hotel selection (\cite{zheng2020using},\cite{DIVSALAR2013150}), clustered POIs (\cite{exposito2019solving},\cite{EXPOSITO2019210}).
\\
\indent In terms of solution methods, meta-heuristic approaches are most commonly used to solve the TTDP and its variant. As claimed in \cite{ruiz2022systematic}, Iterated Local Search (ILS)  or some variations of it (\cite{vansteenwegen2009iterated}, \cite{gavalas2015heuristics},\cite{gavalas2015ecompass}, \cite{doi:10.1287/trsc.1110.0377}) is the most widely applied technique. Indeed, the ILS provides fast and good quality solutions and, therefore, has been embedded in several real-time applications. Other solution methods are: GRASP (\cite{ruiz2021grasp},\cite{EXPOSITO2019210}), large neighboorhod search (\cite{amarouche2020effective}), evolution strategy approach (\cite{karabulut2020evolution}), tabu search (\cite{TANG20051379}), simulated annealing  (\cite{LIN201294}, \cite{LIN2015632}), particle swarm optimization (\cite{DANG2013332}),  ant colony optimisation(\cite{KE2008648}).\\
\indent We finally observe that algorithms solving the TTDP represent one of the main back-end components of expert and intelligent systems designed for supporting tourist decision-making. Among the others they include electronic tourist guides and advanced digital applications such as CT-Planner, eCOMPASS, Scenic Athens, e-Tourism, City Trip Planner, EnoSigTur, TourRec, TripAdvisor, DieToRec, Heracles, TripBuilder, TripSay. A more detailed review of these types of tools can be found in \cite{HAMID2021100337}, \cite{GAVALAS2014319} and \cite{borras2014intelligent}.}

\indent In this paper, we seek to go one step further with respect to the literature by devising insertion and removal operators tailored for a \textit{walk-and-drive} mobility environment. Then we integrate the proposed operators in an iterated local search. A computational campaign on realistic instances show that the proposed approach can handle realistic instances with up to 3643 points of interests in few seconds. The paper is organized as follows.
In section \ref{sec:2} we provide problem definition. In section \ref{sec:3} we describe the structure of the
algorithm used to solve the TTDP. Section \ref{sec:4} and \ref{sec:5} introduce insertion and removal operators to tackle \textcolor{black}{the} TTDP in a \textit{walk-and-drive} mobility environment. Section \ref{sec:6} illustrates how we enhance the proposed approach in order to handle instances with thousands of PoIs. 
In Section \ref{sec:8}, we show the experimental results. Conclusions and further work are discussed in Section \ref{sec:9}.

\section{Problem definition}\label{sec:2}
Let $G=(V,A)$ denote a directed complete multigraph, where each vertex $i\in V$ represents a PoI. Arcs in $A$ are a PoI-based representation of two physical networks: pedestrian network and  road network. \textcolor{black}{Moreover, let $m$ be the length (in days) of the planning horizon.} We denote with $(i,j,mode)\in A$ the connection from PoI $i$ to PoI $j$ with transport $mode\in\{Walk, Drive\}$. Arcs $(i,j,Walk)$ and $(i,j,Drive)$ represent the quickest paths from PoI $i$ to PoI $j$ on the pedestrian network and the road network, respectively. As far as the travel time durations \textcolor{black}{are} concerned, we denote with $t^w_{ij}$ and $t^d_{ij}$  the \textcolor{black}{durations} of the quickest paths from PoI $i$ to PoI $j$ with transport mode equal to $Walk$ and $Drive$, respectively. 
A \textit{score} $P_i$ is assigned to each PoI $i\in V$. Such \textcolor{black}{a} score is determined by taking into account both the popularity of the attraction as well as  preferences of \textcolor{black}{the} tourist. Each PoI $i$ is characterized by a time windows $[O_i,C_i]$ and a visit duration $T_i$.  We denote with $a_i$ the arrival time of the tourist at PoI $i$, with $i\in V$. If the tourist arrives before the opening hour $O_i$, then he/she can wait. \textcolor{black}{Hence,} the PoI visit \textcolor{black}{starts} at time $z_i=max(O_i,a_i)$. The arrival time is feasible if the visit of PoI $i$ can be started before the closing hour $C_i$, i.e. $z_i\leq C_i$. Multiple time windows have been modelled as proposed in \cite{souffriau2013multiconstraint}. Therefore each PoI with more than one time window is replaced by a set of dummy PoI (with the same location and with the same profit) and with one time window each. A \textit{``max-n type"} constraint is added for each set of PoIs to guarantee that at most one PoI per set is visited.

\begin{figure}\label{esempiof}
 \includegraphics[width=150mm ]{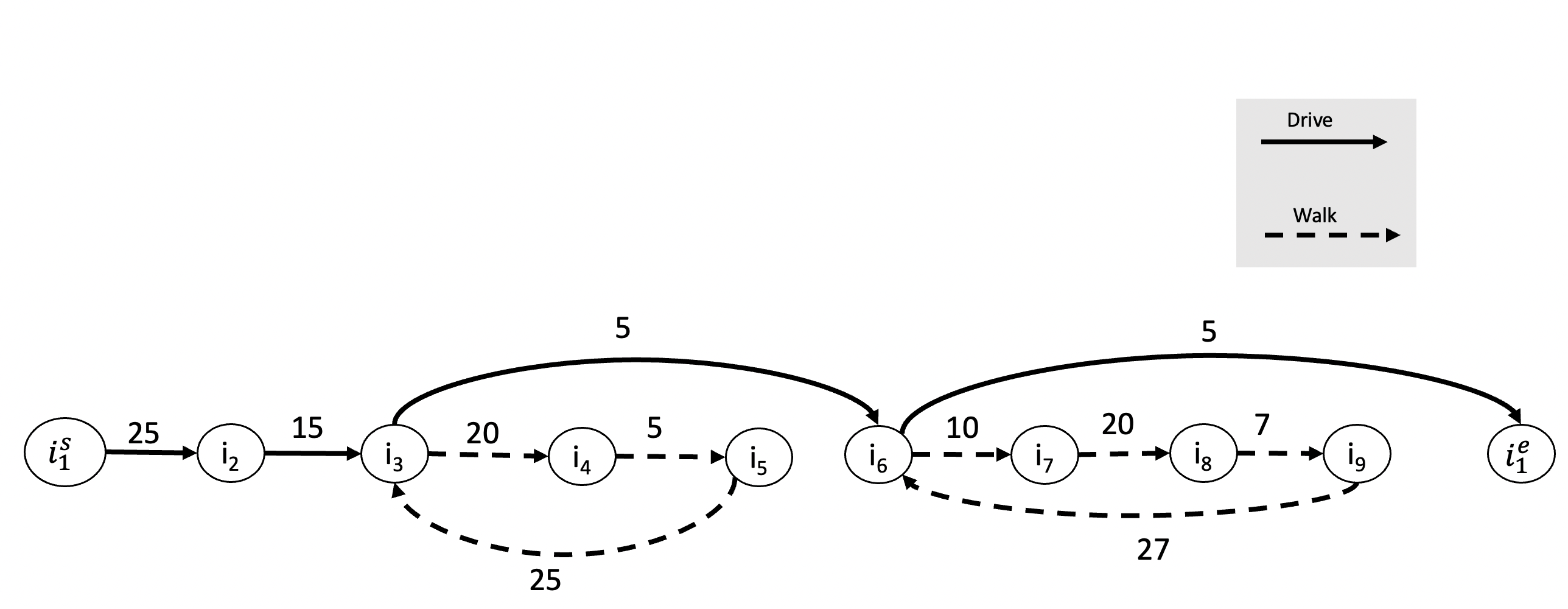}
\caption{Example of a \textcolor{black}{daily itinerary (weights on the arcs indicate travel times)}. }
\end{figure} 

\indent In a \textit{walk-and-drive} mobility environment a TTDP solution consists in the selection of $m$ itineraries, starting and ending to a given initial tourist position. Each itinerary corresponds to a sequence of PoI visits and the transport mode selected for each pair of consecutive PoIs. \textcolor{black}{As an example, Figure 1 depicts the itinerary followed by a tourist on a given day. The tourist drives from node $i^s_1$ to node $i_3$, parks, then follows pedestrian tour $i_3-i_4-i_5$ in order to visit the attractions in nodes $i_3$, $i_4$ and $i_5$. Hence he/she picks up the vehicle parked nearby PoI $i_3$ and drives to vertex $i_6$, parks, then follows pedestrian tour $i_6-i_7-i_8-i_9$ in order to visit the corresponding attractions. Finally the tourist  picks up the vehicle parked nearby PoI $i_6$ and drives to the final destination $i^e_1$ (which may coincide with $i^s_1$).}

Two  parameters model tourist preferences in transport mode selection: $MinDrivingTime$ and $MaxWalkingTime$. Given a pair of PoIs $(i,j)$, we denote with $mode_{ij}$ the  transport mode preferred by the tourist. \textcolor{black}{In the following, we assume that a tourist selects the transportation mode $mode_{ij}$ with the following policy (see Algorithm \ref{alg:alg00}). If} $t^w_{ij}$ is strictly greater than $MaxWalkingTime$, the transport mode preferred by the tourist is $Drive$. Otherwise if $t^d_{ij}$ is not strictly greater than $MinDrivingTime$ (and $t^w_{ij}\leq MaxWalkingTime$), the preferred transport mode is $Walk$. In all remaining cases, the tourist prefers the quickest transport mode. \textcolor{black}{It is worth noting that our approach is not dependent on the mode selection mechanism used by the tourist (i.e., Algorithm \ref{alg:alg00})}. A solution is feasible if the selected PoIs are visited within their time windows and each itinerary duration is not greater than $C_{max}$. The TTDP aims to determine the feasible tour that maximizes the total score of the visited PoIs. Tourist preferences on transport mode selection have been modelled as soft constraints. Therefore, ties on total score are broken by selecting the solution with the minimum number of connections violating tourist preferences.\\
\begin{algorithm}[H]
\SetAlgoLined
  \KwIn{PoI $i$, PoI $j$}
  \KwOut{$mode_{ij} $}
\uIf{$t^w_{ij}>MaxWalkingTime$}{%
  $mode_{ij}\gets Drive$\;
}\uElseIf{$t^d_{ij}\leq MinDrivingTime$}{
   $mode_{ij}\gets Walk$\;
}\Else
{

 \lIfElse{$t^w_{ij}\leq t^d_{ij}$}{$mode_{ij}\gets Walk$}{$mode_{ij}\gets Drive$}
}
  \caption{SelectTransportMode}
  \label{alg:alg00}
\end{algorithm}
\subsection{Modelling transfer}
Transfer connections occur when the tourist switches from the road network to the pedestrian network \textcolor{black}{or} vice versa. Since we assume that \textcolor{black}{tourists always enter} a PoI as a pedestrian, travel time $t^d_{ij}$ has to be increased with transfer times associated to the origin PoI $i$ and the destination PoI $j$. The former models the time required to  pick up the vehicle parked nearby PoI $i$ (\textit{PickUpTime}). The latter models the time required to park and then reach on foot PoI $j$ (\textit{ParkingTime}). During a preprocessing phase we have increased  travel time $t^d_{ij}$ by the (initial) \textit{PickUpTime} and the (final) \textit{ParkingTime}. It is worth noting that a transfer connection also occurs when PoI $i$ is the last PoI visited by a \textit{walking} subtour. In this case, the travel time from PoI $i$ to PoI $j$ corresponds to the duration of a \textit{walk-and-drive} path on the multigraph G: the tourist starts from PoI $i$, reaches on foot the first PoI visited by the walking subtour, then reaches PoI $j$ by driving. In Figure \ref{esempiof} an example of \textit{walk-and-drive} path is $i_5-i_3-i_6$. We observe that the reference application context consists of thousands of daily visitable PoIs. Therefore, it is not an affordable option pre-computing the durations of  $(|V|-2)$ \textit{walk-and-drive} paths associated to each pair of PoIs $(i,j)$. For example in our computation campaign the considered $3643$ PoIs would require more than $180$ GB of memory to store about $5 \cdot 10^{10}$ travel times.  For these reasons we have chosen to reduce significantly the size of the instances by including in the PoI-based graph $G$ only the \textit{PickUpTime} and \textit{ParkingTime}.  As illustrated in the following sections, \textit{walk-and-drive} travel scenarios are handled as a special case of $Drive$ transport mode with travel time computed at run time.

\section{Problem-solving method}\label{sec:3}
Our solution approach is based on  \textcolor{black}{the} \textit{Iterated Local Search (ILS)} \textcolor{black}{proposed in \cite{vansteenwegen2009iterated} for the TOPTW}. To \textcolor{black}{account for a \textit{walk-and-drive} mobility environment, 
we developed a number of extensions and adaptations are discussed in corresponding sections. In our problem, the main decisions amount to determine the sequence of PoIs to be visited and the transport mode for each movement between pairs of consecutive PoIs}. The combination of \textit{walking} subtours and transport mode preferences is the new challenging part of a TTDP defined on a \textit{walk-and-drive} mobility environment. To handle these new features, our ILS contains new contributions compared to the literature. \textcolor{black}{Algorithm \ref{alg:alg0}} reports a general description of ILS. The algorithm is initialized with an empty solution. Then, an improvement phase is carried out by combining a local search and a perturbation step, both described in the following subsections. The algorithm stops when one of the following thresholds is reached:  the maximum number of iterations without improvements or a time limit. 
 The following subsections are devoted to illustrating local search and the perturbation phase.
 

\subsection{Local Search} Given an initial feasible solution (\textit{incumbent}), the idea of \textit{local search} is to explore a neighbourhood of solutions \textit{close} to the incumbent one. Once the best neighboor is found, if it is better than the incumbent, then the incumbent is updated and \textcolor{black}{the} search restarts. 
In our case the local search procedure is an \textit{insertion heuristic}, where the initial incumbent is the empty solution and neighbours are all solutions obtained from the incumbent by adding a single PoI. The neighbourhood is explored in a systematic way by considering all possible insertions in the current solution. As illustrated in Section \ref{sec:4}, the feasibility of neighbour solutions is checked in constant time. As far as the objective function is concerned, we evaluate each insertion as follows. For each  itinerary of the incumbent we consider a (unrouted) PoI $j$, if it can be visited without violating both its time window and the corresponding \textit{max-n type} constraint. Then it is determined the itinerary and the corresponding position with the smallest time consumption. We compute the ratio between the score of the PoI and the \textit {extra time} necessary to reach and visit the new PoI $j$. The ratio aims to model a trade-off between time consumption and score. As discussed in \cite{vansteenwegen2009iterated}, due to time windows the score is considered more relevant than the time consumption during the insertion evaluation. Therefore, the POI $j^*$ with the highest $(score) ^ 2 / (extra$ $time)$ ratio is chosen for insertion. Ties are broken by selecting the insertion with the minimum number of violated soft constraints. After the PoI to be inserted has been selected and it has been determined where to insert it, the affected itinerary needs to be updated as illustrated in Section \ref{sec:5}. This basic iteration of insertion is repeated until it is not possible to insert further PoIs due to the constraint imposed by the maximum duration of the itineraries and by PoI time windows. At this point, we have reached a local optimal solution and we proceed to diversify the search with a  \textit {Solution Perturbation} phase. In Section \ref{sec:6}, we illustrate how we leverage clustering algorithms to identify and explore high density neighbourhood consisting of candidate PoIs with a `good' ratio value.
\subsection{Solution Perturbation}
The perturbation phase has the objective of diversifying the local search, avoiding that the algorithm remains \textit {trapped} in a local optima of  solution landscape. The perturbation procedure aims to remove a set of PoIs occupying consecutive positions in the same itinerary. It is worth noting that the perturbation strategy is adaptive. As discussed in Section \ref{sec:4}, in a multimodal environment a removal might not satisfies the triangle inequality, generating a violation of time windows for PoIs visited later. Since time windows are modelled as hard constraints, the perturbation procedure adapts (in constant time) the starting and ending removal positions so that no time windows are violated. To this aim we relax a soft constraint, i.e. tourist preferences about transport mode connecting remaining PoIs. The perturbation procedure finalizes (Algorithm \ref{alg:alg0} - line \ref{alg:alg0:lab1}) the new solution by decreasing the arrival times to a value as close as possible to the start time of the itinerary, in order to avoid unnecessary waiting times. 
Finally, we observe that the parameter concerning the length of the perturbation ($\rho_d$ in Algorithm \ref{alg:alg0}) is a measure of the degree of search \textit{diversification}. For this reason $\rho_d$ is incremented by 1 for each iteration in which there has not been an improvement of the objective function. If  $\rho_d$  is equal to the length of the longest route, to prevent search from \textit{restarting from the empty solution}, the $\rho_d$ parameter is set equal to 50 $\%$ of the length of the smallest route in terms of number of PoIs. Conversely, if the solution found by the local search is the new \textit{best solution} $s _*$, then  search \textit{intensification} degree is increased and a small perturbation is applied to the current solution  $ s^{\prime} _ {*} $, i.e. $\rho_d$ perturbation is set to 1.\\

\begin{algorithm}[H]
\SetAlgoLined
  \KwData{\textit{MaxIter, TimeLimit}}
  $\sigma_d \gets 1$, $\rho_d \gets1$, $s^{\prime}_{*} \gets\emptyset$, $NumberOfTimesNoImprovement \gets 0$\;
  \While {NumberOfTimesNoImprovement $\leq$ MaxIter Or ElapTime$\leq$ TimeLimit} {
  $s^{\prime}_{*}\gets InsertionProcedure(s^{\prime}_{*})$\;  
\uIf{ $s^{\prime}_{*}$ better than $s_*$}{
$s_* \gets s^{\prime}_{*}$\;
$\rho_d\gets1$\;
$NumberOfTimesNoImprovement \gets 0$\;
}
{$NumberOfTimesNoImprovement \gets NumberOfTimesNoImprovement + 1$\;}
$\rho_d\gets\rho_d + 1$\;
\If {$\rho_d \geq$ \textit{Size of biggest itinerary}}{
 $\rho_d\gets  \max(1,\lfloor$ \textit{(Size of smallest itinerary})$/2\rfloor)$\;
} 
$\sigma_d \gets \sigma_d+ \rho_d$\;
$\sigma_d \gets \sigma_d$ $mod$ \textit{(Size of smallest itinerary)}\;
$s^{\prime}_{*}\gets$\textit{PerturbationProcedure}($s^{\prime}_{*}$,$\rho_d$,$\sigma_d$)\;
\textit{Update ElapTime}\;\label{alg:alg0:lab1}
}
\caption{Iterated Local Search }\label{alg:alg0}
\end{algorithm}

\section{Constant time evaluation framework}\label{sec:4}
\begin{figure}
 \includegraphics[width=150mm ]{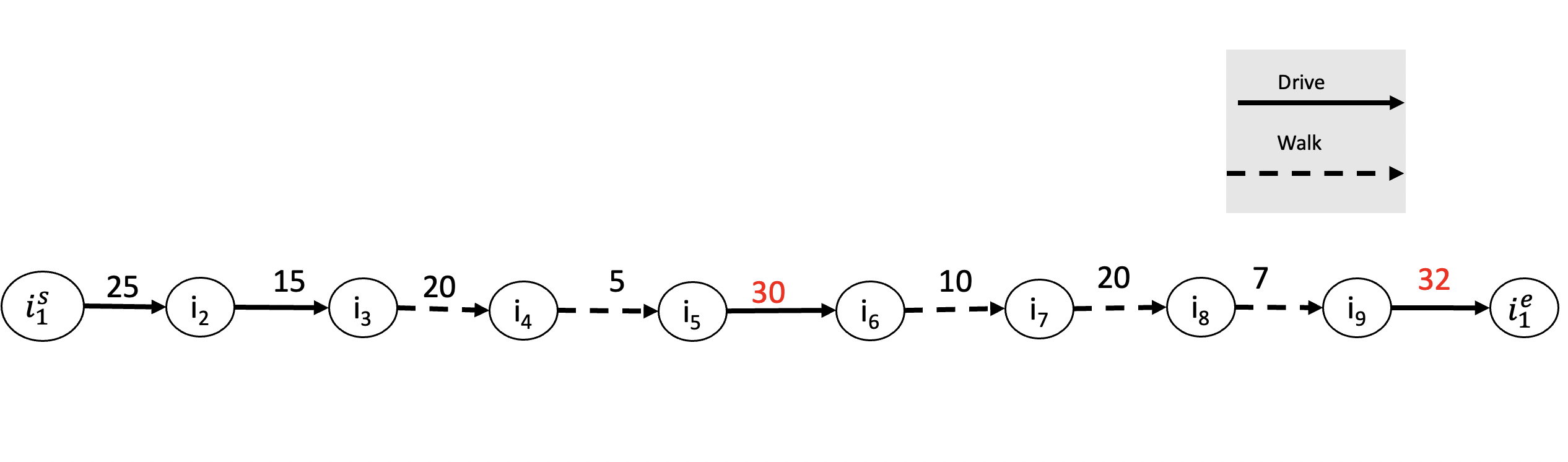}
\caption{Graphical representation of solution encoding of itinerary of Figure \ref{esempiof}. Red travel times refers to duration of \textit{walk-and-drive} paths $(i_5-i_3-i_6)$ and $(i_9-i_6-i^e_1)$. }\label{esempio_enc}
\end{figure}

This section illustrates how to check in constant time the feasibility of a solution chosen in the neighbourhood of $s'_{*}$.
 To this aim the encoding of the current solution has been enriched with additional information.
As illustrated in the following section, such information needs to be updated not in constant time, when 
the incumbent is updated. However this is done much less frequently (once per iteration) than evaluating all solutions in the neighbourhood of the current solution. 
\paragraph{Solution Encoding} We recall that, due to multimodality, a feasible solution has to prescribe for each itinerary a sequence of PoIs and the transport mode between consecutive visits. We encode each itinerary in the solution $s^\prime_*$ as a sequence of PoI visits. Figure \ref{esempio_enc} is a graphical representation of the solution encoding of itinerary of Figure \ref{esempiof}. 
Given two PoIs $i$ and $k$ visited consecutively, we denote with $mode^*_{ik}$ the transport mode prescribed by $s'_*$. We also denote with $t_{ik}$, the travel time needed to move from PoI $i$ to PoI $k$. If the prescribed transport mode is $Drive$, then the travel time $t_{ik}$ has to take properly into account the transfer time needed to switch from the pedestrian network to the road network at PoI $i$. In particular, a transfer connection starting at the origin PoI $i$ might generate a \textit{walking} subtour. For example in the itinerary of Figure \ref{esempiof}, in order to drive from PoI $i_5$ to PoI $i_6$, the tourist has to go on foot from PoI $i_5$ to PoI $i_3$ (\textit{transfer connection}), pick up the vehicle parked nearby PoI $i_3$, drive from PoI $i_3$ to PoI $i_6$ and then park the vehicle nearby PoI $i_6$. In this case we have that $t_{i_5i_6}=t^w_{i_5i_3}+t^d_{i_3i_6}$. To evaluate in constant time the insertion of a new visit between PoIs $i_5$ and $i_6$, we need to encode also subtours. Firstly we maintain two quantities for the $h$-th subtour of an itinerary:  the index of the first PoI and the index of the last PoI denote $FirstPoI_h$ and $LastPoI_h$, respectively.  
For example, the itinerary in Figure \ref{esempiof} has two subtours: the first subtour ($h = 1$) is defined by the PoI sequence $i_3-i_4-i_5$ ($FirstPoI_ {1} = i_3, LastPoI_ {1} = i_5$); the second subtour ($h = 2$) is defined by the PoI sequence $ i_6-i_7-i_8-i_9 $ ($FirstPoI_ {2} = i_6, LastPoI_ {2} = i_9$). We also maintain information for determining in constant time the subtour which a PoI belongs to. In particular, we denote with $S$ a vector of $|V|$ elements: if PoI $i$ belongs to subtour $h$, then $S_i=h$. For the example in Figure \ref{esempiof} we have that $ S_ {i_3} = S_ {i_4} = S_ {i_5} = 1 $, while $ S_ {i_6} = S_ {i_7} = S_ {i_8} = S_ {i_9} = 2$. To model that the remaining PoIs do not belong to any subtour we set $ S_ {i_1} = S_ {i_2} = - 1$. Given two PoIs $i$ and $k$ visited consecutively by solution $s'_*$, the arrival time $a_k$ is determined as follows:
\begin{equation}\label{sol_enc:0}
a_k=z_i+T_i+t_{ik},
\end{equation}
where the travel time $t_{ik}$ is computed by Algorithm \ref{alg:alg2}, according to the prescribed transport $mode$.  If $S_i\neq-1$ and $mode=Drive$, then the input parameter $p$ denote the first PoI of the subtour which PoI $i$ belongs to, i.e. $p=FirstPoI_{S_i}$. If $mode=Walk$ the input parameter $p$ is set to a the deafult value $-1$. Parameter $Check$ is a boolean input, stating if soft constraints are relaxed or not. If $Check$ is $true$, when $mode_{ik}$ violates soft constraints the travel time $t_{ik}$ is set to a large positive value M, making the arrivals at later PoIs infeasible wrt (hard) time-window constraints. In all remaining cases $t_{ik}$ is computed according to the following relationship:
\begin{equation}\label{sol_enc:1}
t_{ik}=t^w+t^d.
\end{equation}
In particular if the prescribed transport mode is \textit{``walk from PoI $i$ to PoI $k$"}, then $t^w=t^w_{ik}$ and $t^d=0$. Otherwise the prescribed transport mode is \textit{``walk from PoI $i$ to PoI $p$, pick-up the vehicle at PoI $p$ and then drive from PoI $p$ to PoI $k$"}, with  $t^w=t^w_{ip}$ and $t^d=t^d_{pk}$.  We abuse notation and when PoI $i$ does not belong to a subtour ($S_i=-1$) and $mode=Drive$, we set  $p=i$ with $t^w_{ii}=0$ and $mode_{ii}=Walk$. A further output of Algorithm \ref{alg:alg2} is the boolean value $Violated$, exploited during PoI insertion/removal to update the number of violated soft constraints.\\
\indent The first six columns of Table \ref{tab:0} report the encoding of the itinerary reported in Figure \ref{esempio_enc}. Tourist position is represented by dummy PoIs $i^s_1$ and $i^e_1$, with a visiting time equal to zero. The arrival time $a_i$ is computed according to (\ref{sol_enc:0}). Column $z_i+T_i$ reports the leaving time with $z_i=max(a_i,O_i)$ and a visiting time $T_i$ equal to 5 time units. All leaving times satisfy time-window constraints, i.e. $z_i\leq C_i$. As far as the timing information associated to the starting and ending PoIs $i^s_1$ and $i^e_1$, they model that the tourist leaves $i^s_1$ at a given time instant (i.e. $a_{i^s_1}=0$), the itinerary duration is 224 time units, with time available for sightseeing equal 320 time units. All connections satisfy soft constraints, since we assume that $MaxWalkingTime$ and $MinDrivingTime$ are equal to 30 and 2 time units, respectively. The last four columns reports details about travel time computations performed by Algorithm \ref{alg:alg2}. Travel time information between PoI $i$ and the next one is reported on the row associated to PoI $i$. Thus this data are not provided for the last (dummy) PoI $i^e_1$. \\
\begin{table}[]\caption{Details of solution encoding for itinerary reported in Figure \ref{esempio_enc}}\label{tab:0}
\resizebox{\textwidth}{!}{%
\begin{tabular}{|c|cc|ccc|cc|cccc|}
\hline
 \multicolumn{6}{|c|}{\textbf{Itinerary}}                                                   & \multicolumn{2}{c|}{\textbf{Time windows}} & \multicolumn{4}{c|}{\textbf{Travel Time Computation}}            \\ \cline{1-12} 
 \textbf{PoI}                 & \textbf{Violated} & $\textbf{mode}^*_{ik}$ & $\textbf{S}_i$ & $\textbf{a}_i$ & $\textbf{z}_i\textbf{+T}_i$ & \textbf{$O_i$}      & \textbf{$C_i$}      & \textbf{p} & $\textbf{t}^w$ & $\textbf{t}^d$ & $\textbf{t}_{ik}$ \\ \hline
$i^s_1$                & False            & Drive           & -1     & 0      & 0          & 0                    & 0                   & $i^s_1$      & 0          & 25         & 25           \\$i_2$                & False              & Drive           & -1     & 25     & 30         & 0                    & 75                  & $i_2$        & 0          & 15         & 15           \\ \hline
$i_3$                & False             & Walk            & 1      & 45     & 55         & 50                   & 115                 & -1           & 20         & 0          & 20           \\
$i_4$                & False             & Walk            & 1      & 75     & 80         & 60                   & 95                  & -1           & 5          & 0          & 5            \\
$i_5$                & False             & Drive           & 1      & 85     & 90         & 60                   & 115                 & $i_3$        & 25         & 5          & 30           \\ \hline
$i_6$                & False             & Walk            & 2      & 120    & 125        & 80                   & 135                 & -1           & 10         & 0          & 10           \\
$i_7$                & False             & Walk            & 2      & 135    & 155        & 150                  & 175                 & -1           & 20         & 0          & 20           \\
$i_8$                & False              & Walk            & 2      & 175    & 180        & 90                   & 245                 & -1           & 7          & 0          & 7            \\
$i_9$                & False             & Drive           & 2      & 187    & 192        & 90                   & 245                 & $i_6$        & 27         & 5          & 32           \\ \hline
$i^e_1$                & -             & -               & -1     & 224    & 224        & 0                    & 320                 & -            & -          & -          & -            \\ \hline
\end{tabular}%
}
\end{table}
 \begin{algorithm}[H]
\SetAlgoLined
\SetKwInput{KwInput}{Input}
\SetKwInput{KwOutput}{Output}
\KwData{M}
  \KwInput{PoI $i$, PoI $k$, $mode$, $Check$, PoI $p$}
  \KwOutput{$t_{ik}$, Violated }
Violated$\gets$ False\;
\uIf{$mode==Walk$}{%
  $t^d\gets 0$\;
  \lIfElse{($Check \wedge mode_{ik}\neq Walk$)} {$t^w\gets M$}{$t^w\gets t^w_{ik}$}
   \lIf{ ($mode_{ik}\neq Walk$)}{Violated$\gets$ True}
} \uElse
{
 \lIfElse{($Check \wedge mode_{ip}\neq Walk$)} {$t^w\gets M$}{$t^w\gets t^w_{ip}$}
 \lIfElse{($Check \wedge mode_{pk}\neq Drive$)} {$t^d\gets M$}{$t^d\gets t^d_{pk}$}
 \lIf{ ($mode_{ip}\neq Walk \vee mode_{pk}\neq Drive$)}{Violated$\gets$ True}

}
$t_{ik}=t^w+t^d$\; 
  \caption{Compute travel time }
  \label{alg:alg2}
\end{algorithm}

\subsection{Feasibility check}
 In describing rules for feasibility checking, we will always consider inserting (unrouted) PoI $j$ between PoI $i$ and $k$. In the following we assume that PoI $j$ satisfies the \textit{max-n type} constraints, modelling multiple time windows. Feasibility check rules are illustrated in the following by distinguishing three main insertion scenarios. The first one is referred to as \textit{basic insertion} and assumes that the extra visit $j$ propagates a change only in terms of arrival times at later PoIs. The second one is referred to as \textit{advanced insertion} and generates a change on later PoIs in terms of both arrival times and (\textit{extra}) transfer time of subtour $S_k\neq-1$. The third one is referred to as a \textit{special case} of the advanced insertion, with PoI $k$ not belonging to any subtour (i.e. $S_k$ is equal to -1). A special case insertion generates a new subtour where PoI $k$ is the last attraction to be visited.\\
  \indent Algorithm \ref{alg:alg3_1} reports the pseudocode of the feasibility check procedure, where the insertion type is determined by $(mode^*_{ik},S_k,mode_{ij},mode_{jk})$. To illustrate the completeness of our feasibility check procedure, we report  in Table \ref{tab:1} all insertion scenarios, discussed in detail in the following subsections. It is worth noting that if $mode^*_{ik}$ is $Walk$ then there exists a \textit{walking} subtour consisting of at least PoIs $i$ and $k$, i.e. $S_k\neq-1$. For this reason we do not detail case 0 in Table \ref{tab:1}.   \\ 
  
%
%
\begin{algorithm}[H]
\SetAlgoLined
  \KwData{PoI $i$, PoI $j$,PoI $k$, incumbent solution $s^*$}
  Compute $Shift_j$ and $Wait_j$\;
\uIf{$mode^*_{ik}=mode_{jk} \wedge (mode_{jk}=Drive \vee mode_{ij}=Walk)$}{
Check Feasibility with (\ref{Feas_ins_1}) and (\ref{Feas_ins_2})\tcp*[f]{Basic Insertion}\;\label{line:1}
}\uElseIf{$S_k\neq-1$}{
		Compute $\Delta_k$ and $Shift_q$\;
		Check feasibility with (\ref{Feas_ins_1_1}), (\ref{Feas_ins_3}) and (\ref{Feas_ins_2})\tcp*[f]{Advanced Insertion}\;\label{line:2}
}\Else{Compute $\Delta_k$ and $Shift_q$\;
	 Check feasibility with (\ref{Feas_ins_3.2.1}), (\ref{Feas_ins_3}) and (\ref{Feas_ins_2})\tcp*[f]{Special Case}\;\label{line:5}
	 }

  \caption{Feasibility check procedure}
  \label{alg:alg3_1}
\end{algorithm}


\begin{table}[]
\caption{Insertion scenarios and their relationships with feasibility check procedures. }\label{tab:1}
\centering
\begin{tabular}{|c|c|c|c|c|}
\hline
Case&$mode^*_{ik}$ & $S_k$    & $(mode_{ij},mode_{jk})$ & Insertion type  \\  \hline


0&  $Walk $ & $=-1$ &-&-\\ \cline{1-5}
\multirow{4}{*}{$1$}& \multirow{4}{*}{$Walk$}& \multirow{4}{*}{$\neq-1$}                 & (Walk,Walk) &Basic  \\  \cline{4-5}
& &                      & (Drive,Drive) &\multirow{3}{*}{Advanced} \\ \cline{4-4}
& &                      & (Walk,Drive)  &                             \\ \cline{4-4}
& &                      & (Drive,Walk)  &                   \\ \hline
\multirow{4}{*}{$2$}&\multirow{4}{*}{$Drive$} & \multirow{4}{*}{$\neq-1$}                 & (Walk,Walk)  &Advanced     \\  \cline{4-5}
& &                      & (Drive,Drive) &  \multirow{2}{*}{Basic}     \\ \cline{4-4}
& &                      & (Walk,Drive) &                              \\ \cline{4-5}
& &                      & (Drive,Walk) & Advanced                \\ \hline
\multirow{5}{*}{$3$}& \multirow{5}{*}{$Drive$} & \multirow{5}{*}{$-1$} & (Walk,Walk)  & Special Case  \\ \cline{4-5}

& &  & (Drive,Drive) & \multirow{2}{*}{Basic}     \\\cline{4-4}
& & & (Walk,Drive)  &                             \\ \cline{4-5}
& & & (Drive,Walk)  & Special Case                \\ \hline
 
\end{tabular}%
\end{table}

\subsubsection{Basic insertion} We observe that in a unimodal mobility environment a PoI insertion is always \textit{basic} \cite{vansteenwegen2009iterated}. 
In a \textit{walk-and-drive} mobility environment an insertion is checked as basic if one of the following conditions hold. 
If PoI $j$ is added to the walking subtour which PoI $i$ and PoI $j$ belong to, i.e. case 1 in Table \ref{tab:1} with $mode_{ij}=mode_{jk}=Walk$. In all other cases we have a basic insertion if it prescribes $Drive$ as transport mode from $j$ to $k$, i.e. case 1 and 2 with $mode_{jk}=Drive$. Five out of 12 scenarios of Table \ref{tab:1} refers to basic insertions. Conditions 
underlying the first three basic insertion scenarios is that $k$ belongs to a walking subtour (i.e. $S_k\neq-1$) and $FirstPoI_{S_k}$ is not updated after the insertion. 
The remaining basic insertions of Table \ref{tab:1} refer to scenarios where before and after the insertion, PoI $k$ does not belong to a subtour. 
All these five scenarios are referred to as basic insertions since the extra visit of PoI $j$ has an impact \textit{only} on the arrival times at later PoIs.
\paragraph{Examples}To ease the discussion, we illustrate two examples of basic insertions for the itinerary of  Figure \ref{esempiof}. Other illustrative examples can be easily derived from Figure \ref{esempiof}.
\begin{itemize}
\item Insert PoI $j$ between  PoI $i=i_3$ and POI $k=i_4$, with $mode_{ij}=Walk$ and $mode_{jk}=Walk$. Before and after the insertion $FirstPoI_{S_k}$ is $i_3$ and, therefore, the insertion has no impact on later transfer connections.
\item Insert PoI $j$ between Insert PoI $i=i^s_1$ and POI $k=i_2$, with $mode_{ij}=Walk$ and $mode_{jk}=Drive$.  Before and after the insertion PoI $i_2$ does not belong to a subtour. Insertion can change only arrival times from PoI $i_2$ on. 
\end{itemize}
\indent To achieve an O(1) complexity for the feasibility check of a basic insertion, we adopt the approach proposed in \cite{vansteenwegen2009iterated} for a unimodal mobility environment and reported in the following for the sake of completeness. We define two quantities for each  PoI $i$ selected by the incumbent solution: $Wait_i$, $MaxShift_i$. We denote with $Wait_i$ the waiting time occurring when the tourist arrives at PoI $i$ before the opening hour:
 $$Wait_i=max\{0, O_i-a_i\}.$$
$MaxShift_i$ represents the maximum increase of start visiting time $z_i$, such that later PoIs can be visited before their closing hour. $MaxShift_i$ is defined by (\ref{ins_1}), where for notational convenience PoI $i+1$ represents the immediate successor of a generic PoI $i$. 
\begin{equation}\label{ins_1}
MaxShift_i=min\{\textcolor{black}{C_i-z_i},Wait_{i+1}+MaxShift_{i+1}\}.
\end{equation}
Table \ref{tab:2} reports values of $Wait$ and $MaxShift$ for the itinerary of Figure \ref{esempiof}.
It is worth noting that the definition of $MaxShift_i$ is a backward recursive formula, initialized with the difference ($C_{max}-z_{max}$), where $z_{max}$ denotes duration of the itinerary.
To check the feasibility of an insertion of PoI $j$ between PoI $i$ and $k$, we compute extra time $Shift_j$ needed to reach and visit PoI $j$, as follows:
\begin{equation}\label{Shift_1}
Shift_j=t_{ij}+ Wait_j+T_j+t_{jk}-t_{ik}. 
\end{equation} 
It is worth noting that travel times are computed by taking into account soft constraints (i.e. input parameter \textit{Check} of Algorithm \ref{alg:alg2} is set equal to \textit{true}).
Feasibility of an insertion is checked in constant time at line \ref{line:1} of Algorithm \ref{alg:alg3_1}  by inequalities (\ref{Feas_ins_1}) and (\ref{Feas_ins_2}).\\
\begin{equation}\label{Feas_ins_1}
Shift_j=t_{ij}+ Wait_j+T_j+t_{jk}-t_{ik}\leq Wait_k+MaxShift_k
\end{equation}
\begin{equation}\label{Feas_ins_2}
z_i+T_i+t_{ij}+ Wait_j\leq C_j.
\end{equation}

\subsubsection{Advanced insertion}

In advanced insertion, the feasibility check has to take into account that the insertion has an impact on later PoIs in terms of both arrival times and  transfer times. Let consider an insertion of a PoI $j$ between PoI $i_2$ and $i_3$ of Figure \ref{esempiof}, with $mode_{i_2j}=mode_{ji_3}=Walk$. The insertion has an impact on the travel time from PoI $i_5$ to PoI $i_6$, i.e. after the insertion travel time $t_{i_5i_6}$ has to be updated to the new value $t^{new}_{i_5i_6}=t^w_{i_5i_2}+t^d_{i_2i_6}$. This implies that we have to handle two distinct feasibility checks. The former has a scope from PoI $i_3$ to $i_5$ and checks the arrival times with respect to $Shift_j$ computed according to (\ref{Shift_1}). The latter concerns PoIs visited after $i_5$ and checks arrival times with respect to $Shift_{i_5}$, computed by taking into account both $Shift_j$ and the new value of $t_{i_5i_6}$. For notational convenience, the first PoI reached by driving after PoI $k$ is referred to as PoI $b$. Similarly, we denote with $q$ the last PoI of the walking subtour, which $k$ belongs to (i.e. if $S_k\neq-1$, then $q=LastPoI_{S_k}$). To check if the type of insertion is advanced, we have to answer the following question: has the insertion an impact on the travel time $t_{qb}$?  To answer it is sufficient to check if after the insertion the value of $FirstPoI_{S_k}$ will be updated, i.e. the insertion changes the first PoI visited by the walking subtour $S_k$.
Five out of 12 scenarios of Table \ref{tab:1} refers to advanced insertions, that is scenarios where $k$ belongs to a \textit{walking} subtour (i.e. $S_k\neq-1$) and $FirstPoI_{S_k}$ is updated after the insertion. 
Algorithm \ref{alg:alg3_1} handles such advanced insertions by checking if one of the following conditions holds.
The insertion of PoI $j$ splits the subtour which PoI $i$ and PoI $j$ belong to, i.e. case 1 in Table \ref{tab:1} with $mode_{ij}=Drive \vee mode_{jk}=Drive$. In all other cases the insertion is checked as advanced if PoI $j$  is \textit{appended} at the beginning of the subtour $S_k$, i.e. case 2 in Table \ref{tab:1} with $mode_{jk}=Walk$. \\
\paragraph{Examples}As we did for basic insertions, we illustrate two advanced insertions for the itinerary of  Figure \ref{esempiof}. Other illustrative examples can be easily derived from Figure \ref{esempiof}.
\begin{itemize}
\item Insert PoI $j$ between  PoI $i=i_7$ and POI $k=i_8$, with $mode_{ij}=Drive$ and $mode_{jk}=Walk$. After the insertion $FirstPoI_{S_k}$ is $j$. Insertion  change $t_{i_9i_1}$ to the new value $t^{new}_{i_9i_1}=t^w_{i_9j}+t^d_{ji_1}$. 
\item Insert PoI $j$ between PoI $i=i_5$ and POI $k=i_6$ , with $mode_{ij}=Walk$ and $mode_{jk}=Walk$.  After the insertion $FirstPoI_{S_k}$ is $i_3$. Insertion change $t_{i_9i_1}$ to the new value $t^{new}_{i_9i_1}=t^w_{i_9i_3}+t^d_{i_3i_1}$. 
\end{itemize}

\indent To evaluate in constant time an advanced insertion, for each PoI $i$ included in solution $s'_*$, three further quantities are defined when $S_k\neq-1$: $\overline{MaxShift}_i$, $\overline{Wait}_i$  and $ME_i$. \\
$\overline{MaxShift}_i$ represents the maximum increase of start visiting time $z_i$, such that later PoIs of subtour $S_i$ can be visited within their time windows. The definition of $\overline{MaxShift}_i$ is computed as follows in (backward) recursive manner starting with $\overline{MaxShift}_q=(C_q-z_q)$.
\begin{equation}\label{over_maxsh}
\overline{MaxShift}_i=min\{\textcolor{black}{C_i-z_i},Wait_{i+1}+\overline{MaxShift}_{i+1}\}.
\end{equation}
$\overline{Wait}_i$ corresponds to the sum of waiting times of later PoIs of subtour $S_i$. We abuse notation by denoting with $i+1$ the direct successor of PoI $i$ and such that $S_{i+1}=S_{i}$. Then we have that 
\begin{equation}\label{over_wait}
\overline{Wait}_i=\overline{Wait}_{i+1}+Wait_{i},
\end{equation}
with $\overline{Wait}_{LastPoI_{S_i}}=Wait_{LastPoI_{S_i}}$.\\
It worth recalling that in a multimodal mobility environment an insertion might propagate to later PoIs a decrease of the arrival times. The maximum decrease that a PoI $i$ can propagate is equal to $\max\{0, a_i-O_i\}$.  $ME_i$ represents the maximum decrease of arrival times that can be propagated from PoI $i$ to $LastPoI_{S_i}$, that is 
\begin {equation}\label{me}
ME_i=\min\{ME_{i+1}, \max\{(0,a_i-O_i)\}\},
\end{equation}
with $ME_{LastPoI_{S_i}}=\max\{(0,a_{LastPoI_{S_i}}-O_{LastPoI_{S_i}})\}$.
 If extra visit of PoI $j$ generates an increase of the arrival times at later PoIs, i.e. $Shift_j\geq0$, then the arrival time of PoI $LastPoI_{S_k}$ is increased by the quantity $max\{0,Shift_j-\overline{Wait}_k\}$. On the other hand if $Shift_j<0$ then the arrival time of PoI $LastPoI_{S_k}$ is decreased by the quantity $\min\{ME_k,|Shift_j|\}$.
Let $\lambda_j$ be a boolean function stating when $Shift_j$ is non-negative:
$$\lambda_j =     \left\{ \begin{array}{rcl}
         1 & & Shift_j\geq0 \\ 
         0  & & Shift_j<0  \\
                \end{array}\right.$$ 
We quantify the impact of extra visit of PoI $j$ on the arrival times of PoI $LastPoI{S_k}$ by computing the value $\Delta_k$ as follows
$$\Delta_{k}=\lambda_j\times\max\{0,Shift_j-\overline{Wait}_{k}\}-(1-\lambda_j)\times\min\{ME_{k},|Shift_j|\}.$$
To check the feasibility of the insertion of PoI $j$ between PoI $i$ and $k$, along with $Shift_j$ we compute $Shift_q$ as the difference between the  new arrival time at PoI $b$ and the old one, that is:
\begin{equation}\label{Sft_q}
Shift_q=t^{new}_{qb}+\Delta_k-t_{qb},
\end{equation}
where $t^{new}_{qb}$ would be the new value of $t_{qb}$ if the algorithm inserted PoI $j$  between PoIs $i$ and $k$.
Feasibility of the insertion of PoI $j$ between PoI $i$ and $k$ is checked in constant time at line \ref{line:2} of Algorithm \ref{alg:alg3_1} by (\ref{Feas_ins_1_1}), (\ref{Feas_ins_3}) and (\ref{Feas_ins_2}).
\begin{equation}\label{Feas_ins_1_1}
Shift_j\leq Wait_k+\overline{MaxShift}_k
\end{equation}
\begin{equation}\label{Feas_ins_3}
Shift_q\leq Wait_b+MaxShift_b.
\end{equation}
Table \ref{tab:2} reports values of $\overline{Wait}$, $\overline{MaxShift}$ and $ME$ for subtours of itinerary of Figure \ref{esempiof}. As we did for basic insertions, travel times are computed by taking into account soft constraints.
 \paragraph{Special case} A special case of the advanced insertion is when PoI $k$ does not belong to a subtour (i.e. $S_k=-1$) in the solution $s'_*$, but it becomes the last PoI of a new subtour after the insertion. Feasibility check rules (\ref{Feas_ins_1_1}) and (\ref{Feas_ins_3}) do not apply since $\overline{MaxShift}_k$, $\overline{Wait}_k$  and $ME_k$ are not defined. In this case,  $\Delta_k$ is computed as follows:
$$\Delta_k= \lambda_j\max(0,Shift_j-Wait_k)-
(1-\lambda_j)\min\{\max\{0,a_k-O_k\},|Shift_j|\}.$$
Then we set $q=k$ and compute $Shift_q$ according to  (\ref{Sft_q}). Feasibility of the insertion of PoI $j$ between PoI $i$ and $k$ is checked in constant time by (\ref{Feas_ins_3.2.1}), (\ref{Feas_ins_3}) and (\ref{Feas_ins_2}).
\begin{equation}\label{Feas_ins_3.2.1}
 Shift_j\leq Wait_k+( C_q-z_q),
\end{equation}

\section{\textcolor{black}{Updating an itinerary}}\label{sec:5}
During the local search after a PoI to be inserted has been selected and it has been decided where to insert the PoI, the affected itinerary needs to be updated. Similarly, during the perturbation phase after a set of selected PoIs has been removed, the affected itineraries need to be updated. The following subsections detail how we update the information maintained to facilitate feasibility checking  when a PoI is inserted  and a sequence of PoI is removed.
\begin{algorithm}[!t]
\caption{Insertion Procedure }
\label{alg:alg4}
\SetAlgoLined
INIT: incumbent solution $s^{\prime}_*$\;
\For {POI j visited by $s^{\prime}_*$}{\label{alg4:2}
 Determine the best feasible insertion with minimum value of $Shift^\prime_j$\;
Compute $Ratio_j$\;
}
Select POI $j^{*}=arg\min\limits_{j}(Ratio_j)$\;\label{alg4:3}
Visit $j^{*}$: Compute $a_{j^{*}}$, $z_{j^*}$, $Wait_{j^{*}}$, $Shift_{j^*}$, $S_{j^*}$\;\label{alg4:5}
Update information of subtours $S_{i^*}$, $S_{k^*}$\; \label{alg4:4}
\lIfElse{Advanced Insertion}{$q^*\gets LastPoI_{S_{k^*}}$, Compute  $Shift_{q^*}$}{$q^*\gets-1$}\label{alg4:1}
$\overline{j}\gets j^{*}$\;
\For (\tcp*[f]{Forward Update}){POI j visited later than $j^{*}$ (Until $Shift_j=0$ $\wedge$ $j \geq q^*$)}{ \label{alg4:6}
 Update $a_j$, $z_j$, $Wait_j$,$S_j$\;
 \lIf{$j\neq q^*$}{Update $Shift_j$}
\lIf{$Shift_j=0$ $\wedge$ $j \geq q^*$}{$\overline{j}\gets j$}\label{alg4:9}
}
\For (\tcp*[f]{Backward Update-Step 1}){POI j visited earlier than $\overline{j}$ (Until $j=FirstPoI_{S_{j^*}}$) }{\label{alg4:7}
Update $MaxShift_j$\;
\lIf{$S_j\neq-1$}{Update  $\overline{Wait}_j$, $\overline{MaxShift}_j$, ${ME}_j$} 
 }
\For (\tcp*[f]{Backward Update-Step 2}){POI j visited earlier than  $FirstPoI_{S_{i^*}}$}{
Update $MaxShift_j$\label{alg4:8}\;
}
Update the number of violated soft constraints\;

\end{algorithm}

\subsection{Insert and Update}
Algorithm \ref{alg:alg4} reports the pseudocode of the proposed insertion procedure. During a major iteration of the local search, we select the best neighbour of the current solution $s'_*$ as follows (Algorithm \ref{alg:alg4} lines \ref{alg4:2}-\ref{alg4:3}). For each (unrouted) PoI $j$ we select the insertion with the minimum value of  $Shift^{\prime}_j=Shift_j+Shift_q$. Then we compute $Ratio_j=(P_j)^2/Shift^{\prime}_j$. The best neighbour is the solution obtained by inserting in  $s'_*$ the PoI $j^*$ with the maximum value of $Ratio_{j^*}$, i.e. $j^*=arg\max\limits_{j}{(P_j)^2/Shift^{\prime}_j}$. Ties are broken by selecting the solution that best fits transport mode preferences, i.e. the insertion with the minimum number of violated soft constraints. The \textit{coordinate} of the best insertion of $j^*$ are denoted with $i^*$, $k^*$. Solution is updated in order to include the visit of $j^*$ (Algorithm \ref{alg:alg4}-lines \ref{alg4:5}-\ref{alg4:4}). If the type of insertion is advanced we determine the value of $Shift_{q^*}$ according to (\ref{Sft_q}) (Algorithm \ref{alg:alg4}-line \ref{alg4:1}). Then, the solution encoding update consists of two consecutive main phases. The first phase  is referred to as \textit{forward update}, since it updates a few information related to visit of PoI $j^*$ and later PoIs.  The \textit{forward update} stops when the propagation of the insertion of $j^*$ has been completely \textit{absorbed} by waiting times of later PoIs (Algorithm \ref{alg:alg4}-lines \ref{alg4:6}-\ref{alg4:9}). 
The second phase is initialized with the PoI $\overline{j}$ satisfying the stopping criterion of the  \textit{forward update}. Such final step is refereed to as \textit{backward update}, since it iterates on PoIs visited earlier than  $\overline{j}$ (Algorithm \ref{alg:alg4}-lines \ref{alg4:7}-\ref{alg4:8}). We finally update the number of violated constraints. As illustrated in the following, new arcs do not violate tourist preferences and therefore after the insertion of $j^*$ the number of violated soft constraints cannot increase.

\paragraph{Solution encoding update}
Once inserted the new visit $j^{*}$ between PoI $i^*$ and PoI $k^*$, we update solution encoding as follows:
\begin{equation}\label{update_1}
a_{j^{*}}=z_i^*+T_i^*+t_{i^*j^{*}}
\end{equation}
\begin{equation}\label{update_2}
Wait_{j^{*}}=\max\{0, O_{j^{*}}-a_{j^{*}}\}
\end{equation}
\begin{equation}\label{update_3}
Shift_{j^{*}}=t_{i^*j^{*}}+ Wait_{j^{*}}+T_{j^{*}}+t_{j^{*}k^*}-t_{i^*k^*}.
\end{equation}
If needed, we update $S_{j^*}$, $FirstPoI_{S_{k^*}}$ and $LastPoI_{S_i{^*}}$. The insertion of $j^*$ propagates a change of the arrival times at later PoIs only if $Shift_{j^*}\neq 0$. We recall that in a multimodal setting, the triangle inequality might not hold. This implies that $j^*$ insertion propagates either an increase (i.e. $Shift_{j^*}>0$) or a decrease (i.e. $Shift_{j^*}<0$)  of the arrival times. Solution encoding of later PoIs is updated according to formula (\ref{update_5})-(\ref{update_7}). For notational convenience we denote with $j$  the current PoI and $j-1$ its immediate predecessor. 
\begin{equation}\label{update_5}
a_{j}=a_{j}+Shift_{{j-1}}
\end{equation}
\begin{equation}\label{update_6}
Shift_{j} =     \left\{ \begin{array}{rcl}
         \max\{0,Shift_{{j-1}}-Wait_{j}\} & & Shift_{j-1}>0 \\ 
         \max\{O_{j}-z_{j},Shift_{j-1}\}  & & Shift_{j-1}<0  \\
                \end{array}\right.
\end{equation}
\begin{equation}\label{update_4}
Wait_{j}=max\{0, O_{j}-a_{j}\}
\end{equation}
\begin{equation}\label{update_7}
z_{j}=z_{j}+Shift_{j}
\end{equation}
At the first iteration, $j$ is initialized with $k^*$ and $Shift_{j-1}=Shift_{j^*}$. In particular (\ref{update_6}) states that after $j$ it is propagated the portion of $Shift_{j-1}$ exceeding $Wait_{j}$, when $Shift_{j-1}>0$. Otherwise $Shift_{j}$ is strictly negative only if no waiting time occurs at PoI $j$ in solution $s'_*$, that is $z_{j}>O_{j}$. If type of insertion is advanced we omit to update $Shift_{q^*}$, since it has been precomputed at line \ref{alg4:1} according to (\ref{Sft_q}).    
The forward updating procedure stops before the end of the itinerary if $Shift_{j}$ is zero, meaning that waiting times have entirely  \textit{absorbed} the initial increase/decrease of arrival times generated by $j^*$ insertion. Then we start the backward update, consisting of two main steps. During the first step the procedure iterates on PoIs visited between the POI  $\overline{j}$,  where the forward update stopped, and $FirstPoI_{S_j^*}$. We update $MaxShift_j$ according to the (\ref{ins_1}) as well as additional information for checking feasibility for advanced insertions. Therefore, if PoI $j$ belongs to a subtour (i.e. $S_j\neq-1$), then we also update $\overline{Wait_j}$, $\overline{MaxShift}_{j}$ and $ME_j$  according to the backward recursive formula (\ref{over_wait}), (\ref{over_maxsh}) and (\ref{me}). The second step iterates on PoI $j$ visited earlier than $FirstPoI_{S_{j^*}}$ and updates only $MaxShift_j$. 
    
\subsection{Remove and Update}
The perturbation procedure aims to remove for each itinerary of the incumbent solution  $\rho_d$ PoIs visited consecutively starting from position $\sigma_d$. Given an itinerary, we denote with $i$ and $k$ respectively the last PoI and the first PoI, that are visited before and after the selected $\rho_d$ PoIs. Let $Shift_{i}$ denotes the variation of total travel time generated by the removal and propagated to PoIs visited later, that is:
$$Shift_{i}=t_{ik}-(a_k-T_i-z_i).$$
In particular when we compute $t_{ik}$ we do not take into account tourist preferences, i.e. in Algorithm \ref{alg:alg2} the input parameter $Check$ is equal to false.
Due to multimodality, the triangle inequality might not be respected by the removal, since it can be propagate either an increase (i.e. $Shift_{i}>0$) or a decrease of the arrival times (i.e. $Shift_{i}<0$). In order to guarantee that after removing the selected PoIs, we obtain an itinerary feasible wrt hard constraints (i.e. time windows), we require that $Shift_{i}\leq0$. To this aims we adjust the starting and the ending removal positions so that it is not allowed to remove portions of  multiple subtours. In particular, if $S_i$ is not equal to $S_k$, then we set the initial and ending removal positions respectively to $FirstPOI_{S_i}$  and the immediate successor of $LastPOI_{S_k}$. In this way we remove subtours $S_i$, $S_k$ along with all the in-between subtours. For example in Figure \ref{esempiof}, if $i$ and $k$ are equal to PoI $i_2$ and $i_4$ respectively, then we adjust $k$ so that the entire first subtour is removed, i.e. we set $k$ equal to $i_6$. Once the selected PoIs have been removed, the solution encoding update steps are the same of a basic insertion. We finally update the number of violated constraints.

\begin{algorithm}[!t]
\caption{Perturbation Procedure  }
\label{alg:alg5}
\SetAlgoLined
INIT: an itinerary of solution $s^{\prime}_*$, i, k\;
$mode=Drive$\;
\uIf{$S_i=S_k$}{
\lIf{$S_i\neq-1$}{$mode\gets Walk$}
} \Else
{
\lIf{$S_i\neq-1$}{$i\gets FirstPoI_{S_i}$}
\lIf{$S_k\neq-1$}{$i\gets$  immediate successor of $LastPoI_{S_k}$}
}
Remove PoIs visited between $i$ and $k$\;
$mode^*_{ik}=mode$\;
$Shift_i\gets t_{ik}-(a_k-z_i-T_i)$\;
Update $a_i$, $z_i$, $Wait_i$\;
\For (\tcp*[f]{Forward Update}){POI j visited later than $i$ (Until $Shift_j=0$)}{ \label{alg5:1}
 Update $a_j$, $z_j$, $Wait_j$\;
\lIf{$Shift_j=0$}{$\overline{j}\gets j$}\label{alg5:2}
}
\For (\tcp*[f]{Backward Update-Step 1}){POI j visited earlier than $\overline{j}$ (Until $j=i$) }{\label{alg5:3}
Update $MaxShift_j$\;
\lIf{$S_j\neq-1$}{Update  $\overline{Wait}_j$, $\overline{MaxShift}_j$, $ME_j$} 
 }
Update $MaxShift_i$\;
\For (\tcp*[f]{Backward Update-Step 2}){POI j visited earlier than  $i$}{
Update $MaxShift_j$\label{alg5:4}
}

\end{algorithm}
\subsection{A numerical example}
We provide a numerical example to illustrate the procedures described so far. We consider the itinerary of Figure \ref{esempiof}. 

\begin{figure}
 \includegraphics[width=150mm ]{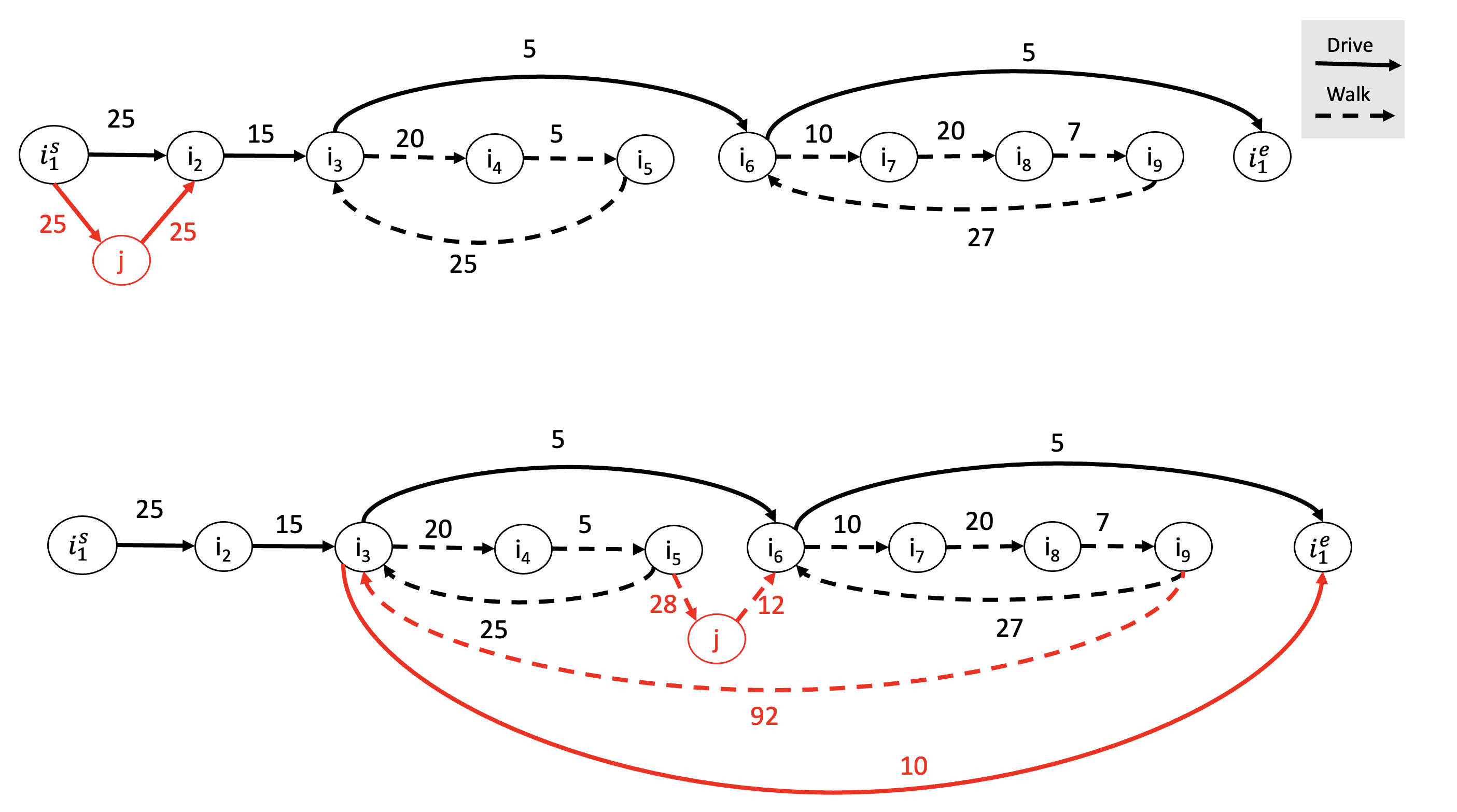}
\caption{Example of infeasible  insertions }\label{Infeas_Ins}
\end{figure}

\begin{figure}
 \includegraphics[width=150mm ]{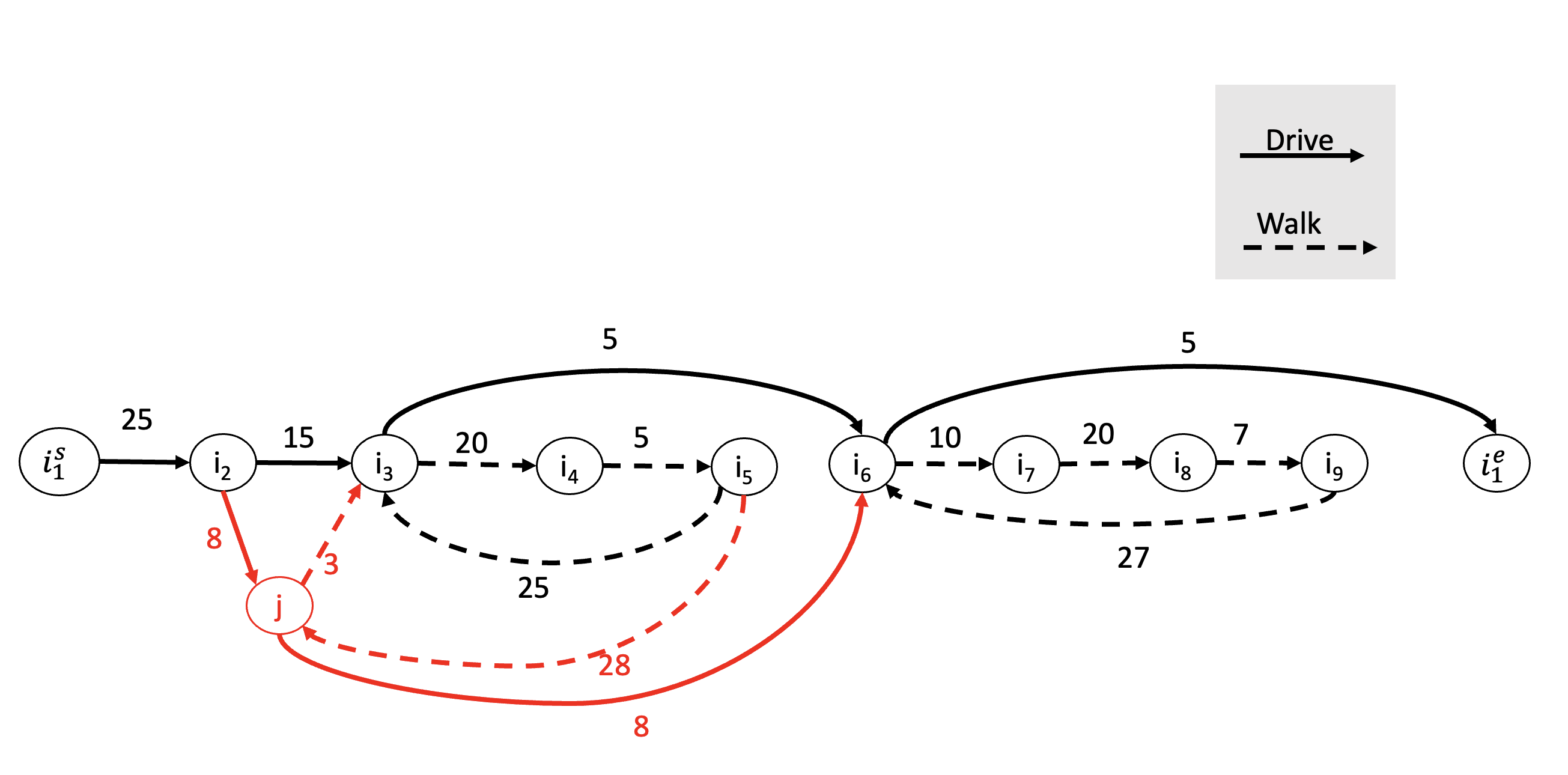}
\caption{Example of feasible  insertion }\label{Feas_Ins}
\end{figure}

In particular we illustrate the feasibility check of the following three insertions for a PoI $j$, with $[O_j,C_j]=[0,300]$ and $T_j=5$.  Durations of arcs involved in the insertion are reported in Figures \ref{Infeas_Ins} and Figure \ref{Feas_Ins}. As reported in Table \ref{tab:2} the itinerary of Figure \ref{esempiof} is feasible with respect to both time windows and soft constraints. As aforementioned, during the feasibility check, all travel times are computed by Algorithm \ref{alg:alg2} with input parameter $Check$ set equal to true.\\
\begin{table}[h!]\caption{Solution encoding with additional information for itinerary of Figure \ref{esempio_enc}}\label{tab:2}
\resizebox{\textwidth}{!}{%
\begin{tabular}{|cccccccc|ccccc|}
\hline
\multicolumn{6}{|c|}{\textbf{Itinerary}}                                                            & \multicolumn{2}{c|}{\textbf{Time Windows}} & \multicolumn{5}{c|}{\textbf{Additional data}}                             \\ \hline
\multicolumn{1}{|c}{\textbf{PoI}} & \textbf{Violated}& $\textbf{mode}^*_{ik}$ & $\textbf{S}_i$ & $\textbf{a}_i$ & \multicolumn{1}{l|}{$\textbf{z}_i\textbf{+T}_i$} & $\textbf{O}_i$               & $\textbf{C}_i$               & $\textbf{Wait}_i$ & $\textbf{MaxShift}_i$ & $\overline{\textbf{Wait}_i}$ & $\overline{\textbf{MaxShift}_i}$ & $\textbf{ME}_i$ \\ \hline
$i_1^s$                &False    & Drive         & -1    & 0     & \multicolumn{1}{l|}{0}         & 0               & 0               & 0    & 0        & -                & -                    & -    \\
$i_2$                &False    & Drive         & -1    & 25    & \multicolumn{1}{l|}{30}        & 0               & 75              & 0    & 20       & -                & -                    & -    \\ \hline
$i_3$                & False   & Walk          & 1     & 45    & \multicolumn{1}{l|}{55}        & 50              & 115             & 5    & 15       & 5                & 20                   & 0    \\
$i_4$                &  False  & Walk          & 1     & 75    & \multicolumn{1}{l|}{80}        & 60              & 95             & 0    & 15       & 0                & 20                   & 15   \\
$i_5$                & False   & Drive         & 1     & 85    & \multicolumn{1}{l|}{90}        & 60              & 115             & 0    & 15       & 0                & 30                   & 25   \\ \hline
$i_6$             & False      & Walk          & 2     & 120   & \multicolumn{1}{l|}{125}       & 80              & 135             & 0    & 15       & 15               & 15                   & 0    \\
$i_7$             &  False     & Walk          & 2     & 135   & \multicolumn{1}{l|}{155}       & 150             & 175             & 15   & 25       & 15              & 25                   & 0    \\
$i_8$              &  False    & Walk          & 2     & 175   & \multicolumn{1}{l|}{180}       & 90              & 245             & 0    & 58       & 0                & 58                   & 85   \\
$i_9$              &  False    & Drive         & 2     & 187   & \multicolumn{1}{l|}{192}       & 90              & 245             & 0    & 58       & 0                & 58                   & 97   \\ \hline
$i_1^e$              & -     & -             & -1    & 224   & \multicolumn{1}{l|}{224}       & 0               & 320             & 0    & 96       & -                & -                    & -    \\ \hline
\end{tabular}
}
\end{table}

\paragraph{Insertion of PoI j between PoI $i^s_1$ and $i_2$ with $mode_{i_1j}=mode_{ji_2}=Drive$} We check feasibility by Algorithm \ref{alg:alg3_1}, with $i=i^s_1$, $k=i_2$. The type of insertion is basic since $mode^*_{ik}= mode_{jk}$ and $mode_{jk}=Drive$.  The feasibility is checked by (\ref{Feas_ins_1}) and (\ref{Feas_ins_2}), that is:
$$Shift_j=t_{ij}+Wait_j+T_j+t_{jk}-t_{ik}=25+0+5+25-25=30\nleq 0+20=Wait_k+MaxShift_k,$$
$$z_i+T_i+t_{ij}+Wait_j=25\leq80=C_j,$$
where travel times $t_{ij}$ and $t_{jk}$ has been computed by Algorithm \ref{alg:alg2} with $p$ set equal to $i^s_1$ and $j$, respectively.  The insertion violates  time window of PoI $i_4$. Such infeasibility is checked through the violation of (\ref{Feas_ins_1}).
\paragraph{Insertion of PoI j between PoI $i_5$ and $i_6$ with $mode_{i_5j}=mode_{ji_5}=Walk$} We check feasibility by Algorithm \ref{alg:alg3_1}, with $i=i_5$, $k=i_6$. The type of insertion is advanced since $mode^*_{ik}\neq mode_{jk}$ and $S_k\neq-1$. We recall that feasibility check consists of two parts. Firstly we check feasibility with respect to (\ref{Feas_ins_1_1}) and (\ref{Feas_ins_2}) that is
$$z_i+T_i+t_{ij}+Wait_j=118\leq300=C_j,$$
$$Shift_j=t_{ij}+Wait_j+T_j+t_{jk}-t_{ik}=15\leq 15=Wait_k+\overline{MaxShift}_k,$$
where travel times have been computed by Algorithm \ref{alg:alg2}, with $p=-1$. 
However the new visit of PoI $j$ is infeasible with respect to soft constraints. As aforementioned this case is encoded as a violation of time windows. Indeed, we compute $Shift_q$ according to (\ref{Sft_q}) with $q=i_9$, $b=i^e_1$, where travel time $t^{new}_{qb}$ is computed by Algorithm \ref{alg:alg2}, with $p=i_3$. Since the tourist has to walk more than 30 time units to pick up the vehicle, i.e. $t^w_{i_9i_3}=92$, then Algorithm \ref{alg:alg2} returns a value $t^{new}_{qb}$ equal to the (big) value M, which violates all time windows of later PoIs.
\paragraph{Insertion of PoI j between PoI $i_2$ and $i_3$ with $mode_{i_2j}=Drive$ and $mode_{ji_3}=Walk$} We check feasibility by Algorithm \ref{alg:alg3_1}, with $i=i_2$, $k=i_3$. The type of insertion is advanced since $mode^*_{ik}\neq mode_{jk}$ and $S_k\neq-1$. The insertion does not violate time windows of PoI $j$ and PoIs belonging to the subtour $S_k$. This is checked by verifying that conditions (\ref{Feas_ins_1_1}) and (\ref{Feas_ins_2}) are satisfied, that is:
$$Shift_j=t_{ij}+Wait_j+T_j+t_{jk}-t_{ik}=1\leq 20=Wait_k+\overline{MaxShift}_k,$$
$$z_i+T_i+t_{ij}+Wait_j=38\leq300=C_j,$$
where $t_{ij}$ and $t_{jk}$ are computed by Algorithm \ref{alg:alg2} with $p=-1$.
Then we check feasibility with respect to closing hours of remaining (routed) PoIs. In particular we compute $Shift_q$ with $q=i_5$, $b=i_6$.  Travel time $t^{new}_{qb}$ is computed with $p=j$. We have that $t^{new}_{qb}=28+8$. Since $Shift_j>0$, then $\Delta_k=\max\{0,Shift_j-\overline{Wait}_k\}=0$. 
$$Shift_q=t^{new}_{qb}+\Delta_k-t_{qb}=36+0-30=6\leq0+15=Wait_b+MaxShift_b.$$
The insertion is feasible since it satisfies also (\ref{Feas_ins_3}).\\
\indent Table \ref{tab:3} shows details of the itinerary after the insertion of PoI $j$ between PoIs $i_2$ and $i_3$. It is worth noting that $Shift_k=0$, but the forward update stops at $\overline{j}=i_7$ since $Shift_q=6$. There is no need to update additional information of later PoIs.\\
\begin{table}[h!]\caption{Details of the itinerary after the insertion}\label{tab:3}
\resizebox{\textwidth}{!}{%
\begin{tabular}{|cccccccc|ccccc|}
\hline
\multicolumn{6}{|c|}{\textbf{Itinerary}}                                                            & \multicolumn{2}{c|}{\textbf{Time Windows}} & \multicolumn{5}{c|}{\textbf{Additional data}}                             \\ \hline
\multicolumn{1}{|c}{\textbf{PoI}} & \textbf{Violated}& $\textbf{mode}^*_{ik}$ & $\textbf{S}_i$ & $\textbf{a}_i$ & \multicolumn{1}{l|}{$\textbf{z}_i\textbf{+T}_i$} & $\textbf{O}_i$               & $\textbf{C}_i$               & $\textbf{Wait}_i$ & $\textbf{MaxShift}_i$ & $\overline{\textbf{Wait}_i}$ & $\overline{\textbf{MaxShift}_i}$ & $\textbf{ME}_i$ \\ \hline
$i_1^s$                &False    & Drive         & -1    & 0     & \multicolumn{1}{l|}{0}         & 0               & 0               & 0    & 0        & -                & -                    & -    \\
$i_2$                &False    & Drive         & -1    & 25    & \multicolumn{1}{l|}{30}        & 0               & 75             & 0    & 13       & -                & -                    & -    \\\hline 
$j$                &False    & Walk         & 1    & 38    & \multicolumn{1}{l|}{43}        & 0               & 300              & 0    & 13       & 4               & 24                    & 0    \\
$i_3$                & False   & Walk          & 1     & 46    & \multicolumn{1}{l|}{55}        & 50              & 115             & 4    & 9       & 4                & 20                   & 0    \\
$i_4$                &  False  & Walk          & 1     & 75    & \multicolumn{1}{l|}{80}        & 60              & 95             & 0    & 9       & 0                & 20                   & 15   \\
$i_5$                & False   & Drive         & 1     & 85    & \multicolumn{1}{l|}{90}        & 60              & 115             & 0    & 9       & 0                & 30                   & 25   \\ \hline
$i_6$             & False      & Walk          & 2     & 126   & \multicolumn{1}{l|}{131}       & 80              & 135             & 0    & 9       & 9               & 9                   & 0    \\
$i_7$             &  False     & Walk          & 2     & 141   & \multicolumn{1}{l|}{155}       & 150             & 175            & 9   & 25       & 9              & 25                   & 0    \\
$i_8$              &  False    & Walk          & 2     & 175   & \multicolumn{1}{l|}{180}       & 90              & 245             & 0    & 58       & 0                & 58                   & 85   \\
$i_9$              &  False    & Drive         & 2     & 187   & \multicolumn{1}{l|}{192}       & 90              & 245             & 0    & 58       & 0                & 58                   & 97   \\ \hline
$i_1^e$              & -     & -             & -1    & 224   & \multicolumn{1}{l|}{224}       & 0               & 320             & 0    & 96       & -                & -                    & -    \\ \hline
\end{tabular}
}
\end{table}

\begin{table}[h!]\caption{Details of the itinerary after the removal}\label{tab:4}
\resizebox{\textwidth}{!}{%
\begin{tabular}{|cccccccc|ccccc|}
\hline
\multicolumn{6}{|c|}{\textbf{Itinerary}}                                                            & \multicolumn{2}{c|}{\textbf{Time Windows}} & \multicolumn{5}{c|}{\textbf{Additional data}}                             \\ \hline
\multicolumn{1}{|c}{\textbf{PoI}} & \textbf{Violated}& $\textbf{mode}^*_{ik}$ & $\textbf{S}_i$ & $\textbf{a}_i$ & \multicolumn{1}{l|}{$\textbf{z}_i\textbf{+T}_i$} & $\textbf{O}_i$               & $\textbf{C}_i$               & $\textbf{Wait}_i$ & $\textbf{MaxShift}_i$ & $\overline{\textbf{Wait}_i}$ & $\overline{\textbf{MaxShift}_i}$ & $\textbf{ME}_i$ \\ \hline
$i_1^s$                &False    & Drive         & -1    & 0     & \multicolumn{1}{l|}{0}         & 0               & 0               & 0    & 0        & -                & -                    & -    \\
$i_2$                &True   & Drive         & -1    & 25    & \multicolumn{1}{l|}{30}        & 0               & 75              & 0    & 50       & -                & -                    & -    \\ \hline
$i_6$             & False      & Walk          & 2     & 32   & \multicolumn{1}{l|}{85}       & 80              & 135             & 48    & 55       & 103               & 55                  & 0    \\
$i_7$             &  False     & Walk          & 2     & 95   & \multicolumn{1}{l|}{155}       & 150             & 175             & 55   & 25       & 55             & 25                   & 0    \\
$i_8$              &  False    & Walk          & 2     & 175   & \multicolumn{1}{l|}{180}       & 90              & 245             & 0    & 58       & 0                & 58                   & 85   \\
$i_9$              &  False    & Drive         & 2     & 187   & \multicolumn{1}{l|}{192}       & 90              & 245             & 0    & 58       & 0                & 58                   & 97   \\ \hline
$i_1^e$              & -     & -             & -1    & 224   & \multicolumn{1}{l|}{224}       & 0               & 320             & 0    & 96       & -                & -                    & -    \\ \hline
\end{tabular}
}
\end{table}

\paragraph{Removal of PoIs between $i_2$ and $i_6$} Table \ref{tab:4} reports details of the itinerary after the removal of PoIs visited between $i_2$ and $i_6$. Travel time $t_{i_2i_6}$ is computed by Algorithm \ref{alg:alg2} with input parameter $Check$ set equal to false. We observe that driving from PoI $i_2$ to PoI $i_6$ violates the soft constraint about $MinDrivingTime$, therefore after the removal the algorithm increases the total number of violated soft constraints.

\section{Lifting ILS performance through \textcolor{black}{unsupervised  learning}}\label{sec:6}
The insertion heuristic explores in a systematic way the neighbourhood of the current solution. Of course, the larger the set $V$ the worse the ILS performance. In order to reduce the size of the neighbourhood explored by the local search, we exploited two mechanisms. Firstly, given the tourist starting position $i^s_1$, we consider an unrouted PoI as candidate for the insertion if it belongs to set:
$$\mathcal{N}_r(i^s_1) = \{i \in V : d(i, i_1^s) \leq r\} \subseteq V$$
where $d: V \times V \rightarrow \mathbb{R}^+$ denotes a non-negative distance function and the radius $r$ is a non negative scalar value. The main idea is that it is likely that the lowest ratio values are associated to PoIs located very far from $i^s_1$. 
We used the Haversine formula to approximate the shortest (orthodromic) distance between two geographical points along the Earth's surface.
The main drawback of this neighbourhood filtering is that a low value of radius $r$ might compromise the degree of diversification during the search. To overcome this drawback we adopt the strategy proposed in \cite{gavalas2013cluster}. It is worth noting that in  \cite{gavalas2013cluster} test instances are defined on an Euclidean space. Since, we use a (more realistic) similarity measure representing the travel time duration of a quickest path, we cannot use k-means algorithm to build a clustering structure. To overcome this limitation we have chosen a hierarchical clustering algorithm.
Therefore, during a preprocessing step we cluster PoIs. The adopted hierarchical clustering approach gives different partitioning depending on the level-of-resolution we are looking at. In particular, we exploited agglomerative clustering which is the most common type of hierarchical clustering. The algorithm starts by considering each observation as a single cluster; then, at each iteration two \emph{similar} clusters are merged to define a new larger cluster until all observations are grouped into a single fat cluster. The result is a tree called dendrogram. The similarity between pair of clusters is established by a linkage criterion: e.g. the maximum distances between all observations of the two sets or the variance of the clusters being merged. 
In this work, the metric used to compute linkage is the walking travel time between pairs of PoIs in the mobility environment: this with the aim of reducing the total driving time.
Given a PoI $i \in V$, we denote with $\mathcal{C}_i$ the cluster label assigned to $i$. $\mathcal{C}_d$ is the cluster containing the tourist starting position.
We enhance the local search so that to ensure that a cluster (different from $\mathcal{C}_d$) is visited at most once in a tour. $\mathcal{C}_d$ can be visited at most twice in a tour: when departing from and when arriving to the depot, respectively.
A PoI $j\in \mathcal{N}_r(i_1^s)$ can be inserted between PoIs $i$ and $k$ in a itinerary $\mathfrak{p}$ only if at least one of the following conditions is satisfied:

\begin{itemize}
\item $\mathcal{C}_i = \mathcal{C}_j \vee \mathcal{C}_k = \mathcal{C}_j$, or
\item $\mathcal{C}_i = \mathcal{C}_k = \mathcal{C}_d \wedge |\mathcal{L}_\mathfrak{p}|=1$, or
\item $\mathcal{C}_i \neq \mathcal{C}_k \wedge \mathcal{C}_j \notin \mathcal{L}_\mathfrak{p} $,
\end{itemize}
where $\mathcal{L}_\mathfrak{p}$ denotes the set of all cluster labels for PoIs belonging to itinerary $\mathfrak{p}$. At first iteration of ILS $\mathcal{L}_\mathfrak{p} = \{\mathcal{C}_d\}$;
subsequently, after each insertion of a PoI $j$, set $\mathcal{L}_\mathfrak{p}$ is enriched with $\mathcal{C}_j$.
In the following section we thoroughly discuss about the remarkable performance improvement obtained, when such cluster based neighbourhood search is applied on (realistic) test instances with thousands of PoIs.

\section{Computational experiments}\label{sec:8}
This section presents the results of the computational experiments conducted to evaluate the performance of our method. We have tested our heuristic algorithm on a set of \textcolor{black}{instances} derived from the pedestrian and road \textcolor{black}{networks} of Apulia (Italy).

All experiments reported in this section \textcolor{black}{were} run on a standalone Linux machine
with an Intel Core i7 processor composed by 4 cores clocked at 2.5 GHz and
equipped with 16 GB of RAM. The machine learning component was implemented in Python (version 3.10). The agglomerative clustering implementations were taken from \textit{scikit-learn} machine learning library. All other algorithms have been coded in Java. \\
Map data \textcolor{black}{were} extracted from OpenStreetMap (OSM) geographic database of the world (publicly available at \url{https://www.openstreetmap.org}). We used the GraphHopper (\url{https://www.graphhopper.com/}) routing engine to precompute all quickest paths between PoI pairs applying an ad-hoc parallel one-to-many Dijkstra for both moving modes (walking and driving). GraphHopper is able to assign a speed for every edge in the graph based on the road type extracted from OSM data for different vehicle profiles: on foot, hike, wheelchair, bike, racing bike, motorcycle and car. 
\textcolor{black}{A fundamental assumption} in our work is that travel times on both driving and pedestrian networks satisfy triangle inequality. In order to satisfy this preliminar requirement, we run \textcolor{black}{the Floyd-Warshall} \cite{floyd1962algorithm,warshall1962theorem} algorithm as a post-processing step to enforce triangle inequality when not met (due to roundings or detours). 
The PoI-based graph consists of $3643$ PoIs.
Walking speed has been fixed to $5$ km/h, while the maximum walking distance is $2.5$ km: i.e. the maximum time that can be traveled on foot is half an hour ($MaxWalkingTime$).
As stated before, we improved the removal and insertion operators of the ILS proposed in order to take into account the extra travel time spent by the tourist to switch from the pedestrian network to the road network. Assuming that the destination has a parking service, we increased the traversal time by car of a customizable constant amount fixed to $10$ minutes ($ParkingTime$). We set the time need to switch from the pedestrian network to the road network equal to at least 5 time minutes ($PickUpTime$). Walking is the preferred mode whenever the traversal time by car is lower than or equal to $6$ minutes ($MinDrivingTime$).
PoI score measures the popularity of attraction. We recall that the research presented in this paper is part of \textcolor{black}{a project} aiming to develop  technologies enabling territorial marketing and tourism in Apulia (Italy). The popularity of PoIs has been extracted from a tourism related Twitter dataset presented in \cite{https://doi.org/10.48550/arxiv.2207.00816}. 
ILS is stopped after $150$ consecutive iterations without improvements or a time limit of one minute is reached. \\
Instances are defined by the following parameters:
\begin{itemize}
\item number of itineraries $m = 1, 2, 3, 4, 5, 6, 7$;
\item starting tourist position (i.e. its latitude and longitude);
\item a radius $r = 10, 20, 50, +\infty$ km for the spherical neighborhood $\mathcal{N}_r(i_1^s)$ around the starting tourist position.
\end{itemize}

\begin{figure}
\centering
 \includegraphics[width=100mm ]{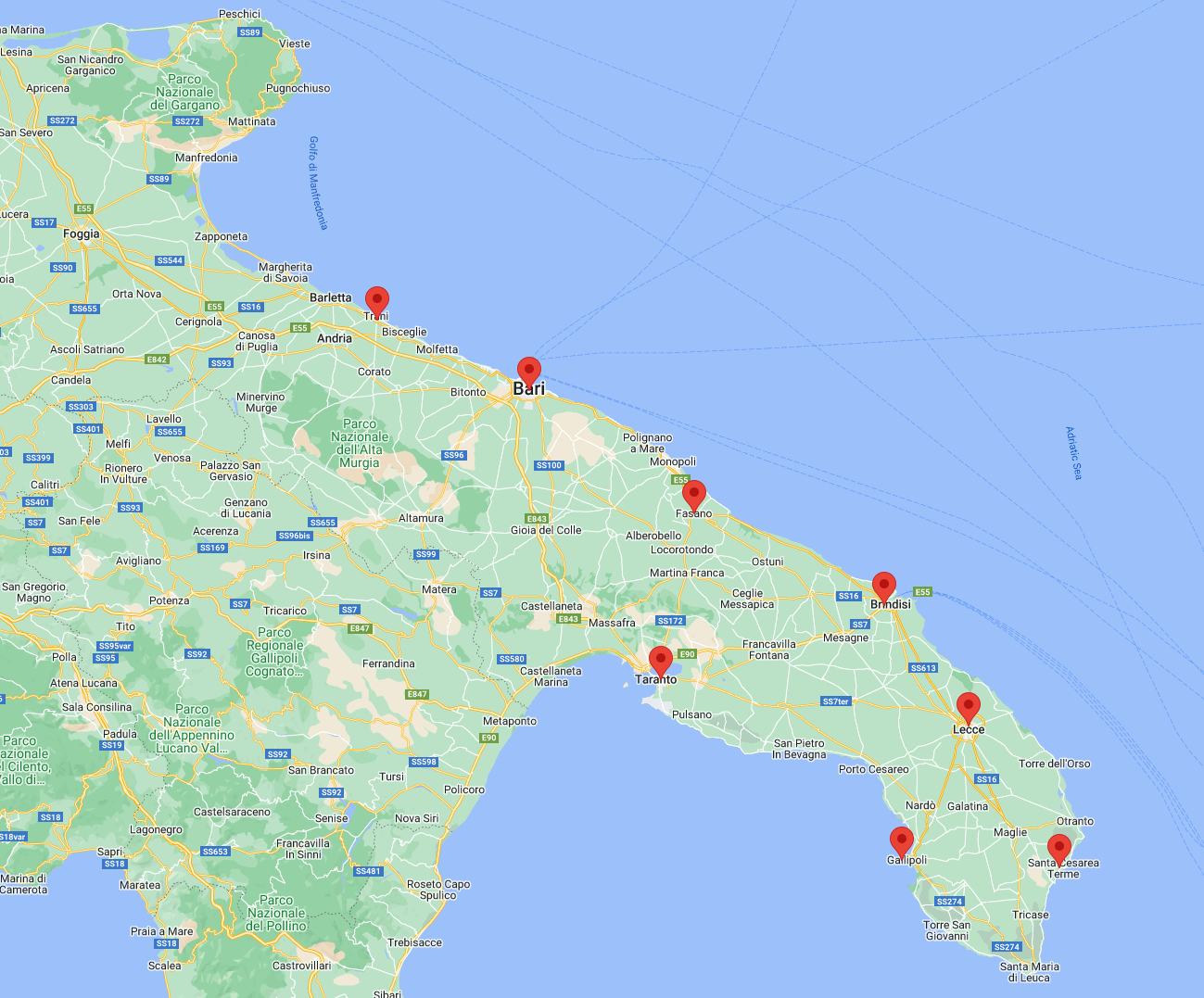}
\caption{Starting positions.\label{fig:apulia}}
\end{figure} 

We considered eight different starting positions along the Apulian territory, as showed in Figure \ref{fig:apulia}. The maximum itinerary duration $C_{max}$ has been fixed to $12$ hours. Every PoI have $0$, $1$ or $2$ opening time-windows depending on current weekday.

\begin{table}[!t]
\normalsize
\caption{Candidate PoIs set size.\label{tableCss}}
\centering
\begin{tabular}{|lc|cccccccc|}
\hline
$r$                       & position & PoIs  & $D_1$   & $D_2$   &$D_3$   &$D_4$   & $D_5$   & $D_6$  & $D_7$   \\ \hline
10                      & 1        & 172  & 257  & 257  & 257  & 203  & 257  & 256  & 148  \\
                        & 2        & 62   & 91   & 91   & 91   & 89   & 91   & 90   & 71   \\
                        & 3        & 172  & 214  & 215  & 176  & 216  & 216  & 136  & 137  \\
                        & 4        & 109  & 118  & 120  & 109  & 120  & 120  & 99   & 99   \\
                        & 5        & 118  & 127  & 132  & 122  & 132  & 132  & 109  & 108  \\
                        & 6        & 79   & 108  & 108  & 107  & 97   & 108  & 106  & 73   \\
                        & 7        & 117  & 140  & 141  & 115  & 141  & 141  & 140  & 140  \\
                        & 8        & 81   & 65   & 82   & 82   & 82   & 82   & 80   & 80   \\ \hline
20                      & 1        & 324  & 507  & 509  & 509  & 385  & 509  & 507  & 254  \\
                        & 2        & 117  & 174  & 172  & 174  & 169  & 174  & 174  & 129  \\
                        & 3        & 301  & 350  & 359  & 312  & 360  & 360  & 264  & 262  \\
                        & 4        & 245  & 266  & 280  & 251  & 280  & 280  & 223  & 222  \\
                        & 5        & 338  & 363  & 390  & 357  & 390  & 390  & 321  & 320  \\
                        & 6        & 262  & 359  & 359  & 346  & 305  & 359  & 354  & 228  \\
                        & 7        & 222  & 260  & 260  & 194  & 261  & 262  & 258  & 253  \\
                        & 8        & 263  & 296  & 329  & 328  & 287  & 329  & 324  & 240  \\ \hline
50                      & 1        & 872  & 1260 & 1289 & 1286 & 1009 & 1289 & 1279 & 712  \\
                        & 2        & 779  & 1010 & 1017 & 928  & 1008 & 1022 & 926  & 776  \\
                        & 3        & 1194 & 1380 & 1437 & 1306 & 1437 & 1441 & 1198 & 1130 \\
                        & 4        & 1267 & 1394 & 1463 & 1311 & 1463 & 1466 & 1202 & 1179 \\
                        & 5        & 1083 & 1185 & 1252 & 1124 & 1254 & 1254 & 994  & 991  \\
                        & 6        & 883  & 1232 & 1230 & 1147 & 1090 & 1235 & 1225 & 832  \\
                        & 7        & 836  & 1083 & 1082 & 938  & 1031 & 1089 & 1081 & 860  \\
                        & 8        & 670  & 875  & 905  & 902  & 768  & 905  & 896  & 606  \\ \hline
$+\infty$ & *        & 3643 & 4591 & 4570 & 4295 & 4521 & 4581 & 4297 & 3781\\
\hline
\end{tabular}
\end{table}

Table \ref{tableCss} summarizes for any radius-position pair:
\begin{itemize}
\item the number of PoIs in the spherical neighborhood $\mathcal{N}_r(i_1^s)$;
\item $D_i$ the number of PoIs opened during day $i$ \textcolor{black}{($i=1, \dots,m=7$, }from Monday to Sunday).
\end{itemize}
When $r$ is set equal to $+\infty$ (last table line) no filter is applied and all $3643$ PoIs in the dataset are candidates for insertion.

Computational results are showed in Table \ref{tableNoClustering}, while Table \ref{tableClustering} reports results obtained with PoI-clustering enabled. Each row represents the average value of the eight instances, with the following headings:
\begin{itemize}
\item DEV: the ratio between total score for the solution and the best known solution;
\item TIME: execution time in seconds;
\item PoIs: number of PoIs;
\item $|S|$: number of walking subtours;
\item SOL: number of improved solutions;
\item IT: total number of iterations;
\item IT$_f$: number of iterations without improvements w.r.t. the incumbent solution;
\item $T^{d}$: total driving time divided by $m \cdot C_{max}$;
\item $T^{w}$: total walking time divided by $m \cdot C_{max}$;
\item $T$: total service time divided by $m \cdot C_{max}$;
\item $W$: total waiting time divided by $m \cdot C_{max}$.
\end{itemize}

Since the territory is characterized by a high density of POIs,  radius $r=50$ km is sufficient to build high-quality tours. 
Furthermore, we notice that the clustering-based ILS greatly improves the execution times of the algorithm, without compromising the quality of the final solution. In particular, the results obtained for increasing $m$ show that, when clustering is enabled, the ILS is able to do many more iterations, thus discovering new solutions and improving the quality of the final solution. 
When  the radius value $r$ is lower than or equal to 50 Km and PoI-clustering is enabled, the algorithm stops mainly due to the iteration limit with $m$ not greater than $5$ itineraries.

The ILS approach is very efficient. The results confirm that the amount of time spent waiting is very small. Itineraries are well-composed with respect to total time spent travelling (without exhausting the tourist). On average, our approach builds itineraries with about $2$ \textit{walking} subtours per day. In particular total walking time and total driving time corresponds  respectively to about 6\% and 20\% of the available time.
On average the visit time corresponds to about the 70\% of the available time. Whilst the waiting time is on average less than 1.5\%.   

We further observe that by increasing the value $r$, the search execution times significantly increase with and without PoI-clustering. With respect to tour quality, clustered ILS is able to improve the degree of diversification on the territory, without remain trapped in high-profit isolated areas.

\begin{table}[!t]
\small
\centering
\caption{Computational results\label{tableNoClustering}}
\begin{tabular}{|ll|ccccccccccc|}
\hline
$m$  &   $r$    & DEV [\%]           & TIME [s]         & PoIs      & $|S|$ & SOL    & IT             & IT$_f$          & $T^{d}$ [\%]   & $T^{w}$ [\%]  & $T$ [\%] & $W$ [\%]    \\\hline
1 & 10    & 16.5          & 0.8           & 18.3          & 1.9             & 2.1          & 155.0          & 150.0          & 13.2          & 6.9          & 78.7          & 1.2          \\
  & 20    & 9.3           & 1.7           & 19.4          & 2.5             & 2.4          & 157.3          & 150.0          & 17.0          & 7.4          & 74.9          & 0.7          \\
  & 50    & 3.4           & 7.0           & 19.8          & 2.4             & 3.6          & 162.5          & 150.0          & 20.9          & 7.2          & 70.9          & 1.0          \\
  & $+\infty$ & 2.7           & 37.5          & 19.5          & 1.3             & 2.6          & 157.1          & 150.0          & 22.8          & 8.2          & 68.1          & 0.9          \\ \hline
2 & 10    & 26.7          & 2.0           & 33.4          & 3.0             & 2.8          & 159.8          & 150.0          & 15.5          & 5.8          & 77.3          & 1.4          \\
  & 20    & 16.8          & 5.3           & 35.0          & 5.1             & 4.5          & 165.8          & 150.0          & 19.5          & 6.0          & 73.2          & 1.3          \\
  & 50    & 5.0           & 27.5          & 38.3          & 4.5             & 8.1          & 183.3          & 150.0          & 20.6          & 6.9          & 71.6          & 0.9          \\
  & $+\infty$ & 1.8           & 60.0          & 38.6          & 4.3             & 5.3          & 89.5           & 74.0           & 24.3          & 6.8          & 67.8          & 1.0          \\\hline
3 & 10    & 31.6          & 3.4           & 46.8          & 5.0             & 4.0          & 176.3          & 150.0          & 14.2          & 5.4          & 78.6          & 1.8          \\
  & 20    & 19.2          & 10.0          & 50.8          & 7.1             & 5.1          & 170.8          & 150.0          & 19.6          & 5.7          & 73.4          & 1.3          \\
  & 50    & 3.2           & 50.5          & 56.0          & 7.8             & 9.6          & 178.3          & 134.3          & 22.3          & 6.5          & 69.9          & 1.4          \\
  & $+\infty$ & 0.7           & 60.0          & 56.5          & 7.9             & 9.0          & 53.6           & 27.5           & 25.6          & 5.6          & 67.8          & 1.0          \\\hline
4 & 10    & 35.0          & 4.9           & 58.9          & 7.8             & 3.1          & 175.1          & 150.0          & 14.5          & 5.0          & 78.8          & 1.6          \\
  & 20    & 21.7          & 16.0          & 65.5          & 9.5             & 7.1          & 190.1          & 150.0          & 20.1          & 5.1          & 73.2          & 1.5          \\
  & 50    & 3.2           & 58.7          & 72.5          & 11.5            & 8.3          & 127.0          & 90.1           & 23.6          & 6.2          & 69.0          & 1.3          \\
  & $+\infty$ & 1.2           & 60.0          & 72.6          & 10.6            & 8.9          & 36.5           & 13.8           & 26.7          & 6.0          & 66.1          & 1.2          \\\hline
5 & 10    & 38.4          & 5.6           & 70.4          & 8.5             & 5.4          & 166.5          & 150.0          & 14.7          & 4.3          & 79.1          & 1.9          \\
  & 20    & 24.4          & 25.4          & 78.8          & 11.3            & 6.3          & 197.5          & 150.0          & 19.6          & 5.0          & 73.7          & 1.7          \\
  & 50    & 2.6           & 60.0          & 89.3          & 13.1            & 7.6          & 88.9           & 59.8           & 23.0          & 6.0          & 69.7          & 1.3          \\
  & $+\infty$ & 1.3           & 60.0          & 89.4          & 12.3            & 6.9          & 26.3           & 11.8           & 24.4          & 5.9          & 68.3          & 1.4          \\\hline
6 & 10    & 41.5          & 7.3           & 80.0          & 9.9             & 4.1          & 191.8          & 150.0          & 14.2          & 4.3          & 78.9          & 2.5          \\
  & 20    & 27.2          & 27.6          & 90.8          & 13.9            & 5.6          & 184.6          & 150.0          & 20.4          & 4.8          & 73.0          & 1.8          \\
  & 50    & 4.6           & 60.0          & 102.6         & 16.0            & 7.3          & 67.3           & 44.9           & 24.4          & 5.7          & 68.6          & 1.3          \\
  & $+\infty$ & 2.0           & 60.0          & 103.9         & 15.5            & 7.6          & 21.8           & 8.4            & 26.7          & 5.8          & 65.9          & 1.5          \\\hline
7 & 10    & 44.1          & 8.0           & 88.1          & 12.9            & 4.5          & 194.4          & 150.0          & 14.1          & 3.9          & 78.4          & 3.6          \\
  & 20    & 28.4          & 34.3          & 104.5         & 15.0            & 6.0          & 180.1          & 150.0          & 19.5          & 4.9          & 73.5          & 2.2          \\
  & 50    & 4.2           & 60.0          & 118.0         & 18.1            & 8.0          & 56.1           & 23.5           & 24.8          & 5.5          & 68.1          & 1.6          \\
  & $+\infty$ & 3.2           & 60.0          & 117.4         & 18.9            & 7.5          & 18.4           & 5.4            & 27.6          & 5.2          & 65.8          & 1.4          \\ \hline
\multicolumn{2}{|c|}{AVG}       & \textbf{15.0} & \textbf{31.2} & \textbf{65.5} & \textbf{9.2}    & \textbf{5.8} & \textbf{133.3} & \textbf{108.7} & \textbf{20.5} & \textbf{5.8} & \textbf{72.2} & \textbf{1.5}\\
 \hline
\end{tabular}
\end{table}

\begin{table}[!t]
\small
\centering
\caption{Computational results with clustering\label{tableClustering}}
\begin{tabular}{|ll|ccccccccccc|}
\hline
$m$  &   $r$    & DEV [\%]           & TIME [s]         & PoIs      & $|S|$ & SOL    & IT             & IT$_f$          & $T^{d}$ [\%]   & $T^{w}$ [\%]  & $T$ [\%] & $W$ [\%]    \\\hline
1 & 10        & 16.7          & 0.6           & 18.3          & 1.9             & 2.1          & 155.9          & 150.0          & 13.7          & 6.5          & 78.6          & 1.2          \\
  & 20        & 9.7           & 0.9           & 19.4          & 2.4             & 2.8          & 157.3          & 150.0          & 16.5          & 6.9          & 75.7          & 0.8          \\
  & 50        & 4.5           & 2.1           & 19.8          & 2.5             & 3.6          & 161.9          & 150.0          & 19.4          & 7.9          & 71.4          & 1.3          \\
  & $+\infty$ & 2.7           & 9.7           & 19.5          & 1.4             & 2.3          & 154.6          & 150.0          & 22.5          & 8.6          & 67.7          & 1.2          \\ \hline
2 & 10        & 26.5          & 1.8           & 33.4          & 4.0             & 4.4          & 176.6          & 150.0          & 15.5          & 5.6          & 77.8          & 1.2          \\
  & 20        & 16.8          & 2.5           & 35.4          & 4.3             & 3.3          & 164.8          & 150.0          & 18.0          & 6.8          & 73.6          & 1.6          \\
  & 50        & 5.3           & 6.6           & 38.4          & 5.1             & 4.6          & 167.8          & 150.0          & 21.5          & 6.7          & 70.5          & 1.3          \\
  & $+\infty$ & 1.5           & 35.5          & 38.9          & 4.3             & 5.4          & 176.5          & 150.0          & 22.0          & 7.6          & 69.5          & 1.0          \\ \hline
3 & 10        & 31.5          & 3.1           & 47.0          & 5.5             & 3.8          & 185.0          & 150.0          & 13.5          & 5.7          & 78.7          & 2.1          \\
  & 20        & 19.2          & 4.9           & 51.3          & 7.1             & 4.4          & 192.0          & 150.0          & 19.2          & 5.8          & 73.5          & 1.5          \\
  & 50        & 3.9           & 14.3          & 56.3          & 8.3             & 9.5          & 183.9          & 150.0          & 22.6          & 6.3          & 70.1          & 1.0          \\
  & $+\infty$ & 1.5           & 59.8          & 56.4          & 6.8             & 7.5          & 156.8          & 111.4          & 23.9          & 6.4          & 68.5          & 1.2          \\ \hline
4 & 10        & 34.7          & 3.9           & 59.4          & 8.3             & 4.1          & 167.9          & 150.0          & 13.9          & 4.9          & 79.4          & 1.8          \\
  & 20        & 22.0          & 7.6           & 64.9          & 9.1             & 5.1          & 191.6          & 150.0          & 19.8          & 5.3          & 73.3          & 1.6          \\
  & 50        & 3.1           & 23.4          & 72.9          & 11.8            & 9.1          & 177.8          & 150.0          & 23.6          & 6.3          & 68.9          & 1.1          \\
  & $+\infty$ & 1.0           & 60.0          & 72.9          & 10.8            & 7.9          & 96.1           & 66.0           & 25.1          & 5.9          & 67.8          & 1.2          \\ \hline
5 & 10        & 38.4          & 7.0           & 70.5          & 9.6             & 4.6          & 211.5          & 150.0          & 14.5          & 5.0          & 78.6          & 1.9          \\
  & 20        & 24.7          & 10.0          & 78.9          & 10.6            & 5.1          & 187.9          & 150.0          & 19.3          & 5.2          & 73.9          & 1.6          \\
  & 50        & 3.2           & 37.2          & 89.1          & 13.9            & 10.3         & 195.3          & 150.0          & 23.2          & 6.3          & 69.1          & 1.4          \\
  & $+\infty$ & 1.3           & 60.0          & 88.4          & 13.4            & 7.9          & 66.6           & 40.9           & 27.2          & 5.4          & 66.4          & 1.1          \\ \hline
6 & 10        & 41.6          & 6.0           & 80.1          & 11.0            & 3.6          & 186.0          & 150.0          & 15.0          & 4.6          & 78.4          & 2.1          \\
  & 20        & 27.1          & 14.0          & 91.5          & 12.6            & 6.5          & 212.8          & 150.0          & 19.9          & 5.1          & 72.9          & 2.1          \\
  & 50        & 3.3           & 49.7          & 104.9         & 16.3            & 10.9         & 187.6          & 131.5          & 24.0          & 5.9          & 68.8          & 1.3          \\
  & $+\infty$ & 1.3           & 60.0          & 105.6         & 17.0            & 10.3         & 51.0           & 29.8           & 26.5          & 5.7          & 66.5          & 1.3          \\ \hline
7 & 10        & 44.2          & 8.3           & 88.0          & 12.6            & 4.1          & 209.6          & 150.0          & 13.8          & 4.4          & 78.0          & 3.8          \\
  & 20        & 28.4          & 14.4          & 104.3         & 16.6            & 4.8          & 169.3          & 150.0          & 19.5          & 4.7          & 73.8          & 2.0          \\
  & 50        & 3.8           & 57.5          & 118.5         & 18.6            & 9.6          & 174.6          & 91.0           & 24.2          & 6.0          & 68.5          & 1.4          \\
  & $+\infty$ & 1.0           & 60.0          & 119.8         & 19.0            & 9.4          & 42.5           & 17.5           & 26.8          & 5.3          & 66.3          & 1.7          \\ \hline
\multicolumn{2}{|c|}{AVG}     & \textbf{15.0} & \textbf{22.2} & \textbf{65.8} & \textbf{9.4}    & \textbf{6.0} & \textbf{162.9} & \textbf{129.9} & \textbf{20.2} & \textbf{6.0} & \textbf{72.4} & \textbf{1.5} \\
 \hline
\end{tabular}
\end{table}

\section {Conclusions}\label{sec:9}
In this paper we have dealt with the tourist trip design problem in a \textit{walk-and-drive} mobility environment, where the tourist moves from one attraction to the following one as a pedestrian or as a driver of a vehicle. Transport mode selection depends on the compromise between travel duration and tourist \textcolor{black}{preferences. We have modelled} the problem as a \textcolor{black}{\textit{Team Orienteering Problem}} with multiple time windows on a multigraph, where tourist preferences on transport modes \textcolor{black}{have been expressed} as soft constraints.The proposed model is novel in the literature. We \textcolor{black}{have also devised} an adapted ILS coupled with an innovative approach to evaluate \textcolor{black}{neighbourhoods} in constant time. To validate \textcolor{black}{our solution approach}, realistic instances with thousands of PoIs \textcolor{black}{have been} tested.  The proposed approach \textcolor{black}{has succeeded} in calculating personalised trips of \textcolor{black}{up to 7 days} in real-time. Future research lines will \textcolor{black}{consider additional aspects, such as traffic congestion and PoI score dependency on visit duration}.


\section*{Acknowledgments}
This research was  supported by Regione Puglia (Italy) (Progetto Ricerca e Sviluppo C\-BAS  CUP B54B170001200007 cod. prog. LA3Z825). This support is gratefully acknowledged.

\bibliography{mybibfile.bib}

\end{document}